%% file: pack4.tex
\documentclass[twoside,openright,11pt,a4paper,epsf]{article}
\usepackage{latexsym}
\usepackage[utf8]{inputenc}
\usepackage{amsmath}
\usepackage{amsthm}
\usepackage{amsfonts}
\usepackage{amssymb}
\usepackage{epsfig}
\usepackage{nicefrac}
\usepackage[all]{xy}
\usepackage{graphicx}

\newcommand{\color}[6]{}
\newcommand{\R}{\mathbb{R}}
\newcommand{\C}{\mathbb{C}}

\renewcommand{\P}{\mathbb{P}}

\newcommand{\Z}{\mathbb{Z}}
\newcommand{\Q}{\mathbb{Q}}

\newcommand{\cc}{\mathbb{\mathcal C}}

\newcommand{\ca}{\mathcal A}
\newcommand{\ce}{\mathcal E}

\newcommand{\cv}{\mathcal V}
\newcommand{\cu}{\mathcal U}

\newcommand{\cn}{\mathcal N}

\newcommand{\cl}{\mathcal{L}}

\newcommand{\nbd}{neighbourhood }
\newcommand{\nbds}{neighbourhoods }
\newcommand{\fonction}[5]
{$$ 
\begin{array}{rcccl}
 #1 & : & #2 & \longrightarrow &#3 \\
    &   & #4 & \longmapsto &#5 
\end{array}
$$}

\newcommand{\res}{\textnormal{Res}\,}
\newcommand{\reg}{\textnormal{\tiny Reg}}

\newcommand{\rond}[1]{\overset{\vspace*{-1pt}\circ}{#1}}
\newcommand{\priv}{\backslash}
\newcommand{\lra}{\longrightarrow}
\newcommand{\hra}{\hookrightarrow}

\newcommand{\om}{\omega}
\newcommand{\eps}{\varepsilon}
\renewcommand{\phi}{\varphi}

\newcommand{\st}{\textnormal{\small st}}
\newcommand{\vol}{\text{Vol}\,}

\newcommand{\wdt}[1]{\widetilde{#1}}
\newcommand{\cqfd}{\hfill $\square$ \vspace{0.1cm}\\ }
\newcommand{\sbull}{{\tiny $\bullet$ }}
\newcommand{\ds}{\displaystyle}

\renewcommand{\bar}[1]{\overline{#1}}
\newcommand{\pd}{\textnormal{\small PD}}

\newtheorem{definition}{Definition}[section]
\newtheorem{thm}{Theorem}
\newtheorem*{thm*}{Theorem}
\newtheorem{prop}[definition]{Proposition}
\newtheorem{lemma}[definition]{Lemma}
\newtheorem*{question}{Question}
\newtheorem{cor}[definition]{Corollary}
\newtheorem{rk}[definition]{Remark}

\title{ \vspace*{0cm}Symplectic packings in dimension $4$ and singular curves.}
\author{Emmanuel Opshtein. \footnote{Partially supported by ANR project "Floer Power" ANR-08-BLAN-0291-03.}}
\date{}
\begin{document}
\maketitle
\begin{abstract}
The main goal of this paper is to give constructive proofs of several existence results for symplectic embeddings. The strong relation 
between symplectic packings and singular symplectic curves, which can be derived from McDuff's inflations on the blow-ups, is revisited 
through a new inflation technique that lives at the level of the manifold. As an application, we explain constructions of maximal 
symplectic packings of $\P^2$ by $6$, $7$ or $8$ balls.
\end{abstract}
\section{Introduction.}
In \cite{gromov}, M. Gromov shows a striking relation between curves and symplectic embeddings. Namely,  the symplectic curves   may give  obstructions to symplectic embeddings beyond the volume 
constraint. For packings ({\it i.e.} embeddings of disjoint balls), the basic idea is to produce symplectic curves 
through the centers of the balls of the packing, within a prescribed homology class. Such a curve automatically 
gives an obstruction to the size of the balls. In \cite{mcpo},
McDuff and Polterovich proved a converse statement for less than nine balls in $\P^2$ : in particular, the only symplectic 
obstructions are given by such curves. This result hints at an even deeper relation between packings and symplectic curves, which has been mostly confirmed by McDuff {\it via} her blow-up and inflation techniques \cite{mcduff1,lamc2,mcduff3}. For instance, the following theorem follows from these techniques (see also \cite{biran4}) :
\begin{thm}\label{packball1}
Let $(M^4,\om)$ be a symplectic manifold with rational class ($\om\in H^2(M,\Q)$). Then, there is a symplectic packing 
by closed disjoint balls 
$$
\coprod_{i=1}^p \overline{B^4(a_i)}\overset{\om}{\hra} M \hspace{,5cm} (a_i\in \Q)
$$ 
if and only if there exists an irreducible symplectic curve $\Sigma$ Poincaré-Dual to $k[\om]$, with $p$ nodal
singularities of multiplicities $(ka_1+1,\dots,ka_p+1)$. 
\end{thm}
These "$+1$" are only here to conform to the tradition that a theorem should be true rather than almost true and they can mostly be disregarded
 in a first approximation (see also theorem \ref{packball}).
The rough idea is that symplectic blow-up and blow-down establish a correspondence between symplectic forms on blow-ups and 
existence of packings, while the inflation technique - based on curves of the blow-up, hence on singular curves of $M$ - produces symplectic forms on the blow-up. This strategy works well
for theoretic result, but not for constructions. The main lack for effectiveness lies in the above mentioned 
correspondence, which is far from being explicit. Several works have therefore been concerned with explicit constructions 
of symplectic packings, with completely independant techniques coming from toric geometry \cite{schlenk,traynor}. On this side of the story,  the "obstruction curves"  for these packings - the symplectic curves responsible for the obstruction - are completely discarded. In short, inflation provides existence results from the obstruction curves but fail to be effective, while these toric constructions hide the symplectic aspects of the packings, since they do not describe the 
obstruction curves at all. The main object of this paper is to reprove theorem \ref{packball1} in an effective way : 
given a symplectic curve with prescribed singularity, we construct the desired symplectic packing. The idea is to use 
Liouville vector fields in order to "inflate" directly into the manifold rather than on the blow-up. This more direct 
approach leads to a more accurate version of theorem \ref{packball1} (see section \ref{statementsection})
and allows for instance to construct maximal symplectic packings of $\P^2$ by six, seven and eight balls (up to five balls, the constructions are already available, see \cite{moi3}).

\begin{thm}\label{explicitballs}
There exist symplectic packings of $\P^2$ by six, seven and eight open balls of capacity $\nicefrac{2}{5}$, $\nicefrac{3}{8}$ and $\nicefrac{6}{17}$ respectively (the maximal capacities). These open balls each intersect the predicted obstruction curve along a finite number of Hopf discs. 
\end{thm}
Of course, we will not only prove this theorem but also explain the constructions of these maximal packings. The boundary regularity of these packings will be also quickly adressed. 

Close to the packing problem, the question whether an ellipsoid embeds into a symplectic manifold has also been 
largely considered. The question is really about flexibility of symplectic embeddings : how much can such a map 
fold an euclidean piece such as an ellipsoid \cite{schlenk2,mcsc} ? Again, McDuff used the inflation process, this time in successive blow-ups
to get the optimal results \cite{mcduff4} : in a certain class of symplectic manifold (containing $\P^2$ or the $4$-ball), 
embedding an ellipsoid is equivalent to some specific packing problem. Again however, this answer is of 
theoretic nature, and does not provide constructions, even from ball packings. The next result is in the wake of theorem \ref{packball1} and 
mostly follows from \cite{mcduff4}. 

\begin{thm}\label{packell} 
Let $(M^4,\om)$ be  a symplectic manifold with rational class. Then, the closed ellipsoid 
$$
\overline{E(a,b)}:=\tau\overline{E(p,q)}\hspace{,5cm} (p,q\in \Z,\; \tau\in \Q,\; \gcd(p,q) =1)
$$
symplectically embeds into $M$ if and only if there exists an irreducible symplectic curve $\Sigma$, Poincaré-Dual 
to $k[\om]$, with an ordinary singular point, with local symplectic model 
$$
\prod_{i=1}^{k\tau} (z_1^q-\alpha_j z_2^p),\hspace{,5cm} (\alpha_j\in \C).
$$
\end{thm}
The existence of an ellipsoid embedding also depends on the existence 
of a singular symplectic curve, now with a multi-cusp. The number of branches and the singularity type of the cusps
are respectively responsible for the size  and the shape of the ellipsoid. Again, our aim is to give an effective proof. It will also be clear that ellipsoid packings can  be considered, and simply need curves with several singularities.

The paper is organized as follows. First, we give more precise versions of theorem \ref{packball1}, which will be usefull in practice, for instance 
to prove theorem \ref{explicitballs}. Then, we explain the main idea of the paper by sketching a proof for theorem \ref{packball1}. In section \ref{liouvillesec}, we explain some properties and constructions of Liouville forms that will be needed to perform our inflations. Section \ref{ballsec} and \ref{ellsec} are devoted to the proofs of the theorems on balls packings and ellipsoid embeddings respectively, letting 
aside technical assertions about Donaldson's method. We deal with this last point in section \ref{donsec} and end up the paper with some remarks and questions.

\paragraph{Notations :}We adopt the following conventions throughout this paper : \vspace{-,1cm}
\begin{itemize}
\item[-] {\it All} angles will take value in $\R/ \Z$. In other terms an angle $1$ is a full turn in the plane, and the integral of the form $d\theta$ over a circle around the origin in the plane is $1$.\vspace{-,3cm}
\item[-] The standard symplectic form on $\C^2=\R^{4}$ is $\om_\st:=\sum dr_i^2\wedge d\theta_i$, where $(r_i,\theta_i)$ are polar coordinates on the plane factors. With this convention, the euclidean ball of radius $1$ has capacity $1$. \vspace{-,3cm}
\item[-] A Liouville form $\lambda$ of a symplectic structure $\om$ is a one-form satisfying $\om=-d\lambda$.
The standard Liouville form on the plane is $\lambda_\st:=-r^2d\theta$. \vspace{-,3cm}
\item[-] A symplectic ball or ellipsoid is the image of a Euclidean ball or ellipsoid in $\C^n$ by a symplectic embedding. \vspace{-,3cm}
\item[-] The Hopf discs of a Euclidean ball in $\C^2$ are its intersections with complex lines. \vspace{-,3cm}
\item[-] $\ce(a,b)$ denotes the $4$-dimensional ellipsoid $\{a^{-1}|z|^2+b^{-1}|w|^2<1\}\subset \C^2(z,w)$. Because of our 
normalizations, its Gromov's width is $\min(a,b)$. \vspace{-,7cm}
\item[-] When there is no ambiguity for the symplectic form $\om$ on a manifold $M$, we often abbreviate $(M,\tau \om)$
by $\tau M$.
\end{itemize}

\paragraph{Aknowledgements :} I wish to thank V. Kharlamov and D. Auroux for the usefull discussions we had.
\section{More accurate statements.}\label{statementsection}
In this paragraph, we present three variations on theorem \ref{packball1}, where we either release an hypothesis (such as the irreducibility) or get a slightly more accurate result (constructing open objects of maximal size for instance). These
precision may not be of fundamental importance, but they will be usefull in practice. The first remark is that irreducibility 
can be released for cheap. In the theorem below, an $n$-cross of size $a$ refers to a figure composed of $n$ symplectic 
discs of area $a$ which all intersect transversally at exactly one point. A symplectic model is simply the intersection of $B^4(a)$ with $n$ complex lines in $\C^2$ (with the same K\"ahler angles), which we  call a {standard} $n$-cross. 
\begin{thm}\label{packball} Let $(M^4,\om)$ be a rational symplectic manifold. Then, 
\begin{itemize}
\item[i)] If there exists a symplectic packing by closed balls
$$
\overline{B(a_1)}\sqcup\dots\sqcup \overline{B(a_p)}\hra M, \hspace{1cm} a_i\in \Q,
$$
then for all sufficiently large $k$ for which the $ka_i$ are integers, there exists a symplectic irreducible curve Poincaré-Dual to $k[\om]$ in $M$, with exactly $p$ nodes of multiplicity $(ka_1,\dots,ka_p)$.  
\item[ii)] Given a symplectic curve $\Sigma$ Poincaré-Dual to $k\om$ with $p$ nodes of degree $(a_1k_1,\dots,a_pk_p)$, one can construct a symplectic packing by open balls
$$
B(a_1)\sqcup\dots \sqcup B(a_p) \hra M 
$$
provided that  $k_i\geq k$ for each $i$ and that the curve $\Sigma$ contains $p$ disjoint closed $k_ia_i$-crosses of size $a_i$ around each node. When $\Sigma$ is irreducible and $k_i=k$, this last condition is equivalent to the volume obstruction.
\end{itemize}
\end{thm}
Theorem \ref{packball} obviously implies theorem \ref{packball1}. Notice that when only the interior of the crosses are disjoint, theorem \ref{packball} ii) allows to construct balls of size $a_i-\eps$ for arbitrarily small $\eps. $ The next variation is concerned with this construction 
part only, and states that we can even take $\eps=0$. 
This apparently innocent precision is rather expensive for the proof : it needs the introduction of the notion of tame Liouville form and some analysis (not very difficult however). This is the reason why we state it separately : the proof will be given independently, so that the reader can discard this technical part at first.  The main purpose of this precision 
is to allow for constructions of maximal packings by open balls of maximal size, and not only by balls of approximately optimal size. In my opinion, the importance of this point is of conceptual nature. The impossibility of reaching the 
limit size would mean that a deep rigidity phenomenon appears, which would need an explanation - recall that 
prior to the present paper, the only available proof for the {\it existence} of maximal packings by open balls in $\P^2$
relied on a deep result on symplectic isotopies \cite{mcduff3}. Theorem \ref{packopen} ensures that such a 
phenomenon does not happen in general, and that nothing deep hides here. 

\begin{thm}\label{packopen}
Let $(M,\om)$ be a rational symplectic manifold. From a symplectic curve $\Sigma$ Poincaré-Dual to $k[\om]$
with $p$ nodes of multiplicities $(k_1a_1,\dots,k_pa_p)$, one can construct a symplectic packing by $p$ open
balls  of capacity $a_i$, provided that $k_i\geq k$ for all $i$ and that $\Sigma$ contains $p$ disjoint {\bf open}
$k_ia_i$-crosses of size $a_i$ around each node.
\end{thm}
The last variation we give concerns the rationality hypothesis. The main classical importance of this hypothesis 
is the existence of symplectic polarizations in the sense of Biran \cite{biran3} (curves Poincaré-Dual to a multiple of the symplectic form). Of course, 
no such curve exists if the symplectic form is irrational, but a singular notion of polarization was 
defined in \cite{moi4}.  
\begin{definition}
A polarization of a symplectic manifold $(M,\om)$ is a union of weighted symplectic curves 
${\bf \Sigma}:= \{(\Sigma_l,\tau_l),\; l=1\dots n\}$ ($\tau_l\in \R^+$), such that 
$$
\sum_{l=1}^n \tau_l\pd(\Sigma_l)=[\om],
$$
the intersections between different curves are positive, and the singularities of each curve has the symplectic type 
of a complex singularity (these last two conditions can be replaced by assuming $J$-holomorphicity for $\Sigma$  
but almost complex curves play no role here). The curve $\Sigma:=\cup \Sigma_i$ will be called the total curve of the polarization.
\end{definition}
Although the theorems of this paper are easier to state and feel in the rational setting, this definition makes the rationality 
hypothesis completely irrelevant. 

\begin{thm}\label{packballirr}
Let $(M,\om)$ be a symplectic $4$-manifold. Then, there is a symplectic packing by closed disjoint balls
$$
\coprod_{i=1}^p \overline{B^4(a_i)}\overset{\om}{\hra} M 
$$ 
if and only if there exists a polarization $(\Sigma_l,\tau_l)_{l=1\dots n}$ of $M$ with the following properties :
\begin{itemize}
\item[\sbull] Its total curve $\Sigma$ has $p$ nodes $(x_1,\dots,x_p)$ of multiplicities $k^i$ ($i\in [1,p]$). We also denote by 
$k_l^i$ the number of branches of $\Sigma_l$ through $x_i$ (so $\sum k_l^i=k^i$) 
\item[\sbull] $\Sigma$ possesses $p$ disjoint closed $k^i$-crosses of size $a_i$ centered at $x_i$,
\item[\sbull] For each $i$, $\sum k_l^i\tau_l> a_i$. 
\end{itemize}
\end{thm}
It will be clear from the proof that the same precision on construction of open balls from open crosses will also hold in this situation. 
Moreover, all these variations also hold for the problem of ellipsoid embeddings.

\section{Sketch of the proof.}
Before going into proofs, let us illustrate the basic idea of the paper by explaining (not proving) theorem 
\ref{packball1} for one ball.

\paragraph{i)} In one direction, the argument is the following. Assume you have a ball $B$ of capacity $a$ in $M$. 
We wish to prove that for $k$ sufficiently large, there is a curve, Poincaré-Dual to $k\om$ whose intersection with $B$ is 
exactly $k$ Hopf discs. Now by Donaldson's result, there are curves Poincaré-Dual to $k\om$ for $k$ large. The observation is 
that the work is almost done. These curves naturally pass through $\partial B$ along $ka$ Hopf circles, because we know (from 
McDuff blow-up process) that we can think of $\partial B$ as a  curve of symplectic area $a$, whose points correspond to the
Hopf circles of $\partial B$. Now cutting the part of this curve inside the ball, and pasting the $ka$ corresponding Hopf discs gives the 
desired singular curve. 

\paragraph{ii)} Conversely, given an irreducible curve $\Sigma$ Poincaré-Dual to $k\om$ with exactly one singularity, composed 
of $ka$ transverse branches,  we need to construct a ball of capacity $a$. The construction mainly proceeds from the three following observations.
\begin{itemize}
\item[\sbull] First, provided that $a$ satisfies the volume constraint $\nicefrac{a^2}{2}<\vol M$, we can find a $ka$-cross $X$ of size $a$ in $\Sigma$ ($ka$ discs in $\Sigma$ of area $a$, which intersect one another exactly at the singular point). Indeed, 
$$
\ca_\om(\Sigma)=\int_\Sigma \om \overset{\text{PD}}{=} \int_M \om\wedge k\om=2k\vol M>ka^2=ka\cdot a.
$$
By standard \nbd theorems, there exists an embedding $\phi$ of a \nbd of the cross $X_{ka}:=\cup\Delta_i:=\cup \{z=\alpha_i w\}\cap B^4(a)$ into a \nbd of the cross $X$ in $M$. The desired ball will be obtained by inflating this standard \nbd of the cross in $M$ through well-chosen Liouville vector fields.
\item[\sbull] Next, there is a (contracting) Liouville vector field on $B^4(a)\priv X_{ka}$, which points outwards the $\Delta_i$
and whose flow, restricted to any \nbd of $X_{ka}$, but considered for all positive time (when it is defined), recovers the whole 
of $B^4(a)$.
\item[\sbull] Finally, there is a Liouville vector field on $M\priv \Sigma$, also pointing outwards $\Sigma$. Since $M$ is assumed to be closed, this vector field is defined for all positive time.
\end{itemize}
The construction of the embedding then goes as follows. Take any point in $B$. From this point, follow negatively the flow of the 
Liouville vector field, until you reach the domain of definition of the local embedding $\phi$. Note $\tau$ the amount of time you had 
to flew. Then use $\phi$ to send the point inside $M$, and flow positively along the Liouville vector field in $M$ for a time $\tau$. 
Provided this map is well-defined, it is defined on $B^4(a)$ by the second point above, and it is symplectic. For this map to be 
well-defined, we need that the Liouville vector field in $M$ is an extension of the push-forward of the Liouville vector field in 
$B^4(a)$ by $\phi$. This is a cohomological constraint which involves the residue of the corresponding Liouville forms, and which 
is easy to deal with.

\section{Generalities on Liouville forms.}\label{liouvillesec}
As we explained in the previous paragraph, the central objects of this paper are Liouville forms, on the manifolds and on the objects 
we wish to embed. The aim of the present section is to present the main features and examples of these Liouville forms which we will 
need subsequently.
\subsection{Residues of a Liouville form.}
Consider a symplectic manifold $(M,\om)$ equiped with a decomposition 
$$
[\om]=\sum_{i=1\dots n} \tau_i \pd(\Sigma_i),
$$
where $\tau_i$ are real numbers and $\Sigma_i$ are possibly singular symplectic curves. Note that the symplectic form is exact on $M\priv \Sigma$, so there are Liouville forms on $M\priv \Sigma$. 
Fix also a smooth, local disc fibration above the regular part $\Sigma_i^\reg$ of $\Sigma_i$, with the additional property that each disc must intersect $\Sigma_i$ orthogonaly with respect to the symplectic form. We mostly view these fibrations as families 
$D_i(p)$ of symplectic discs, smoothly parameterized by $p\in \Sigma_i^\reg$. In each such disc, we choose a collection of loops
$\gamma_i^\eps(p)$ which approach $p$ as $\eps$ goes to $0$. 
\begin{lemma} \label{extliouville}
Let $\cv$ be a tubular \nbd of $\Sigma$ in $M$, and $\beta$ a Liouville form on $\cv\priv \Sigma$. Then :
\begin{itemize}
\item[i)] If $p\in \Sigma_i$, the numbers 
$$
\beta(\gamma_i^\eps(p)):=\int_{\gamma_i^\eps(p)}\beta
$$
have a limit $b_i(p)$ when $\eps$ goes to zero.
\item[ii)] The numbers $b_i(p)$  depend neither on $p\in \Sigma_i$, nor on the chosen fibration. The number $b_i$ will be 
called  the residue of $\beta$ at $\Sigma_i$.
\item[iii)] If $\beta$ has residues  $(\tau_i)$ at $\Sigma$, then the form $\beta$ extends to a Liouville form on $M\priv \Sigma$.
\end{itemize}
Similarly, if $\cv$ is a regular \nbd in $M$ of any simply connected open set $U$ of $\Sigma$, 
and  $\beta$ is a Liouville form on $\cv\priv \Sigma$, then :
 \begin{itemize}
\item[iv)]  The restriction of $\beta$ to any proper open subset of $\cv$  
extends to a Liouville form on $M\priv \Sigma$, as long as $\beta$ has residues $(\tau_i)$ at $\Sigma$.
\end{itemize}
\end{lemma}
Here are some illustrations of this lemma.
\begin{enumerate}
\item If $\Sigma$ is an irreducible polarization in $M$, any Liouville form $\beta$ on $M\priv \Sigma$ verifies
$$
\lim_{\eps\to 0} \beta(\gamma_i^\eps(p))=\tau, \hspace{,5cm}\text{ where }[\om]=\tau\pd(\Sigma).  
$$ 
For instance, the residue of a Liouville form on $\P^2$ minus an irreducible curve of degree $k$ is $1/k$.
\item By contrast, when there is a linear dependence between the Poincaré-Dual classes of the $\Sigma_i$, a lot of residues are allowed. 
If $\Sigma_i$ are (possibly singular) curves of $\P^2$ of degree $d_i$, there is a Liouville form on $M\priv \cup\Sigma_i$ with residues $\nicefrac{p_i}{d_i}$ at $\Sigma_i$ as soon as $\sum p_i=1$ (simply average the Liouville forms on $M\priv \Sigma_i$). In particular, there is a Liouville form on $\P^2\priv \cup \Sigma_i$ whose residues at any branch of $\Sigma$ is $\nicefrac{1}{d_1+\dots+d_n}$. 
\end{enumerate}

\noindent {\it Proof of lemma \ref{extliouville} :} The point i) is obvious : $\beta(\gamma_i^\eps(p))$ converges to 
$\ca_\om(D_i(p))-\beta(\partial D_i(p))$. In order to prove ii), observe that the symplectic orthogonals to $D(p)$ induce a 
connection on our local fibration above $\Sigma_i^\reg$. For two points $p,q\in \Sigma_i^\reg$, consider a path $\alpha:[0,1]\to \Sigma_i^\reg$ between $p$ and $q$, and parallel transport  a circle $C_p\subset D_i(p)$ along $\alpha$ (it is possible provided $C_p$ is close enough to $p$). Call $S$ the surface obtained by gluing the three following pieces :
\begin{itemize}
\item[\sbull] the cylinder $\ds A:=\underset{t\in [0,1]}{\large \cup} P_\alpha^t(C_p)$,  where $P^t_\alpha$ denote the parallel transport,
\item[\sbull] the disc $D_1$ bounded by $C_p$ in $D_i(p)$,
\item[\sbull] the disc $D_2$ bounded by $C_q:=P^1(C_p)$ in $D_i(q)$.
\end{itemize}   
Now by definition of our connexion, $A$ is a lagrangian cylinder. This point has two consequences. On one hand, since $S$ 
is null-homologic, the symplectic area of $D_1$ and $D_2$ are the same. On the other hand, since $A$ lies in $\cv\priv \Sigma$, on which $\beta$ is a well-defined Liouville form for $\om$, we get 
$$
\beta(C_p)=\beta(C_q)+\ca_\om(A)=\beta(C_q).
$$ 
We therefore see that $b_i(p)=b_i(q)$, since 
$$
b_i(p)=\beta(C_p)-\ca_\om(D_1)=\beta(C_q)-\ca_\om(D_2)=b_i(q).
$$
In order to see that the residue does not depend on the chosen fibration, consider two fibrations $D^1_i(p),D_i^2(p)$, construct a third one which contains the discs $D^1_i(p)$ and $D^2_i(q)$ and apply the last result to this new fibration. 

In order to prove iii), consider a Liouville form $\lambda$ on $M\priv \Sigma$. By ii), this form has residues $\tau_i+f_i$ 
at $\Sigma_i$. When the residues of $\lambda$ are  the $\tau_i$, the extension of $\beta$ is obvious. Indeed,
$\beta -\lambda$ is closed on $\cv\priv\Sigma$ and has no period by hypothesis, so it is even exact. Writing $\beta-\lambda=dh$, 
any extension $\wdt h$ of $h$ to $M\priv \Sigma$ provides an extension $\wdt \beta:=\lambda+d\wdt h$ of $\beta$. Assume now 
that the $f_i$ do not vanish. The first observation is that 
\begin{equation}\label{vanpd}
\sum f_i\pd(\Sigma_i)=0.
\end{equation}
Indeed, if $C$ is any surface in $M$,
$$
\begin{array}{ll}
\ds \int_C\om &= \ds \sum_{i=1}^n \res (\lambda,\Sigma_i)C\cdot \Sigma_i=\sum_{i=1}^n (\tau_i+f_i)C\cdot \Sigma_i\\
 &  \ds = \sum_{i=1}^n \tau_i\pd(\Sigma_i)\cdot [C]=\sum_{i=1}^n \tau_i C\cdot\Sigma_i, 
\end{array}\\
$$
so $\sum f_i\pd(\Sigma_i)\cdot C=0$.

Consider now a collection of closed $2$-forms $\sigma_i$ representing the class $\pd(\Sigma_i)$. They are obviously exact on 
$M\priv \Sigma_i$, so there exist $1$-forms on $M\priv \Sigma_i$ such that $\sigma_i=-d\alpha_i$. It is easy to check that $\alpha_i$
also has a well-defined residue at $\Sigma_i$, and even that $\res(\alpha_i,\Sigma_i)=1$. Since by (\ref{vanpd}), the form $\sum f_i\sigma_i$ is exact on $M$, we can write 
$$
d(\sum f_i\alpha_i)=d\theta,
$$
where $\theta$ is a $1$-form defined on the whole of $M$. The form $\alpha:=\sum f_i\alpha_i-\theta$ is therefore closed on $M\priv \Sigma$, with residues $f_i$ at $\Sigma_i$. The form $\lambda-\alpha$ is thus a Liouville form on $M\priv \Sigma$, with residues 
$\tau_i$ on $\Sigma_i$, and we are in the previous situation. The proof of iv) is completely similar. \cqfd

\subsection{Liouville forms on  balls.}\label{liouvilleballsec}
The relevant Liouville forms on the balls are restrictions of global Liouville forms on $\C^2$ obtained in the following way. 
The symplectic form on $\C^2$ is $\om_\st=dR_1\wedge d\theta_1+dR_2\wedge d\theta_2$, where
$(R_1:=|z_1|^2,\theta_1,R_2:=|z_2|^2,\theta_2)$ are polar coordinates on $\C^2(z_1,z_2)$, so for  $d_0:=\{z_2=0\}\subset \C^2$,
the $1$-form
$$
\lambda_{\tau}:=(\tau-R_2)d\theta_2-R_1d\theta_1
$$ 
is Liouville, has residue $\tau$ at $d_0$. Any convex combination of pull backs of such Liouville forms by unitary maps gives a Liouville 
form on the complement of lines in $\C^2$ with identified residues. 
\begin{prop}\label{liouvilleball}
Let $X:=\Delta_1\cup\dots\cup \Delta_n$ be a standard $n$-cross in $B(a)\subset \C^2$. If $\tau_1,\dots,\tau_n\in \R^+$ are such that $\sum \tau_i\geq a$,  there exists a Liouville form $\lambda$ on $B(a)\priv X$, with residues 
$$
\res(\lambda,\Delta_i)=\tau_i,
$$
and whose associated Liouville vector field has the following properties :
\begin{itemize}
\item[i)] It is not defined on $X$, but it points outwards the cross.
\item[ii)] Its flow is radius increasing : $X_\lambda\cdot r>0$. In particular, the negative flow of any point in $B(a)$ is well-defined until it reaches one of the $\Delta_i$. In other terms, the basin of repulsion of  $X$, defined by  
$$
\left\{p\in \rond B(a)\,|\, \exists \tau\in \R^+ \text{ with }
\left|\begin{array}{l}
\Phi_{X_\lambda}^{-t}(p) \text{ exists for } t\in [0,\tau[\\
\ds \lim_{t\to \tau} \Phi^{-t}_{X_\lambda}(p)\in X
\end{array}\right.
  \right\} 
$$ 
is exactly $\rond B(a)$. 
\item[iii)] More basically, if $\cu$ is any \nbd of $X$ in $B(a)$, 
$$
\bigcup_{t\geq 0}\Phi^t_{X_\lambda}(\cu)=\rond B(a).
$$
\end{itemize}
\end{prop}
\noindent {\it Proof : } Observe that when $\kappa\geq a$, the $1$-form $\lambda_\kappa$ defined above is associated to the 
vector field 
$$
X_\kappa:=(\kappa-R_2)\frac{\partial}{\partial R_2}-R_1\frac{\partial}{\partial R_1}, 
$$
with $X_\kappa\cdot (R_1+R_2)=\kappa-R_1-R_2$, so $\Phi^t_{X_0}$ increases the radius inside $B(a)$ (it is tangent to $\partial B(a)$ exactly when $\kappa=a$). Moreover, since 
$$
\frac{\partial}{\partial R_2}=\frac{1}{2r_2}\frac{\partial}{\partial r_2},
$$
$X_\kappa$ explodes and points outwards near $d_0$. Given now $(\Delta_i),(\tau_i)$ as in proposition \ref{liouvilleball}, consider 
unitary transformations $u_i$ taking $d_0$ to $d_i:=\langle \Delta_i\rangle$, positive weights $\mu_i$ (with $\sum \mu_i=1$), 
and real numbers $\kappa_i\geq a$ such that $\mu_i\kappa_i=\tau_i$.   
Since the $u_i$ are symplectic and radius-preserving,  the Liouville  form 
$$
\lambda:=\sum \mu_i u_{i*}\lambda_{\kappa_i}
$$
has the properties announced.\cqfd

\subsection{Liouville forms on ellipsoids.}
Ellipsoids can be presented as ramified symplectic coverings of balls. Lifting the  Liouville 
forms from the balls, we now produce the analogue Liouville forms in ellipsoids.
Fix an ellipsoid $E(a,b)$ with $\nicefrac{a}{b}\in \Q$, and put $a:=\tau p$, $b:=\tau q$, where $p,q$ are  relatively prime integers. 
\begin{lemma}\label{covellball} The map 
\fonction{\Phi}{B(\tau pq)}{E(a,b)=\tau E(p,q)}{(R_1,\theta_1,R_2,\theta_2)}{\big(\frac{R_1}{q},q\theta_1,\frac{R_2}{p},p\theta_2\big)}
is a symplectic covering of degree $pq$, ramified over the coordinate axis $\{R_1=0\}$ and $\{R_2=0\}$. It is invariant by the 
diagonal action of $\Z_p\times \Z_q\approx \big(e^{\nicefrac{2i\pi k}{p}},e^{\nicefrac{2i\pi l}{q}}\big)$ on $\C^2$. 
\end{lemma}
\begin{definition}\label{defcone}
Let $\alpha:=|\alpha | e^{i\phi}$ be a complex number. We denote by $\Delta_\alpha\subset E(a,b)$ the cone defined 
by 
$$
\left\{
\begin{array}{l}
R_1=|\alpha| R_2\\
q\theta_1=p\theta_2+\phi
\end{array} .  \right.
$$
It is a symplectic surface, smooth except at the origin (the vertex of the cone), and whose intersection with 
$\partial E(a,b)$ is a characteristic leaf.
\end{definition}
A straightforward computation shows that the preimage of $\Delta_\alpha$ is a union of $pq$ lines, invariant 
by the action of $\Z_p\times \Z_q$. Averaging a Liouville form with residue $\kappa$ around one of this disc (given by proposition \ref{liouvilleball}) by this action, we therefore get a Liouville form on $B(\tau pq)$, which descends to a Liouville form on $E(a,b)\priv \Delta_\alpha$, with the same nice features as for the Liouville forms on the ball. For instance, if the residue at $\Delta_\alpha$ is at least $\tau$,  then the residues of the Liouville forms on $B(\tau pq)$ are at least $\tau$, and the Liouville vector field is radius increasing. The Liouville vector field on the ellipsoid therefore  increases the function $qR_1+pR_2$, hence the function $\nicefrac{R_1}{a}+\nicefrac{R_2}{b}$.  Averaging the different Liouville forms associated to different $\Delta_\alpha$, we get the analogue of proposition 
\ref{liouvilleball} for ellipsoids. 
 
\begin{prop}\label{liouvillell} Let $\alpha_1,\dots,\alpha_n$ be $n$ complex numbers and  
$\Delta_{\alpha_1},\dots,\Delta_{\alpha_n}\subset \tau E(p,q)$ their associated cones. If $\tau_1,\dots,\tau_n\in \R^+$ are such that $\sum \tau_i\geq \tau$,  there exists a Liouville form $\lambda$ on $\tau E(a,b)\priv \Delta_{\alpha_1}\cup\dots\cup\Delta_{\alpha_n}$, with residues 
$$
\res(\lambda,\Delta_{\alpha_i})=\tau_i,
$$
and whose associated Liouville vector field has the following properties :
\begin{itemize}
\item[i)] It is not defined on $\cup\Delta_{\alpha_i}$, but it points outwards these curves. 
\item[ii)] Its flow is "radius increasing". Namely, if $R:=\frac{R_1}{a}+\frac{R_2}{b}$, $X_\lambda\cdot R>0$. In particular, the negative flow of any point in $\tau E(a,b)$ is well-defined until it reaches one of the $\Delta_{\alpha_i}$. In other terms, the basin of repulsion of  $\cup \Delta_{\alpha_i}$, defined by  
$$
\left\{p\in  \tau \rond E(a,b)\,|\, \exists T\in \R^+ \text{ with }
\left|\begin{array}{l}
\Phi_{X_\lambda}^{-t}(p) \text{ exists for } t\in [0,T[\\
\ds \lim_{t\to T} \Phi^{-t}_{X_\lambda}(p)\in \Delta_{\alpha_i}
\end{array}\right.
  \right\} 
$$ 
is exactly $\tau \rond E(a,b)$. 
\item[iii)] More basically, if $\cu$ is any \nbd of $\cup \Delta_{\alpha_i}$ in $\tau E(a,b)$, 
$$
\bigcup_{t\geq 0}\Phi^t_{X_\lambda}(\cu)=\tau \rond E(a,b).
$$
\end{itemize}
\end{prop}

\subsection{Tameness.}
In this paragraph, we define a regularity notion for Liouville forms, which will be central in producing maximal open packings or 
embeddings. 

\begin{definition}[Angular form] Let $\Sigma$ be a codimension-$2$ submanifold in $M$ and $\cn_\Sigma$ its normal bundle in $M$.
Endow $\cn_\Sigma$ with a hermitian metric and a hermitian connexion $\alpha$. An angular form 
 around $\Sigma$ is the push-forward of $\alpha$ by any smooth local embedding of $(\cn_\Sigma,0_{\cn_\Sigma})$ into 
 $(M,\Sigma)$. More informally, it is a $1$-form which is locally the $d\theta$ of local polar coordinates around $\Sigma$. 
 
 For symplectic submanifold, we further impose a compatibility condition: $\alpha$ must be 
 positive on the small closed loops that turn positively around $\Sigma$.
\end{definition} 

\begin{definition}[Tameness] Let $(M^4,\Sigma)$ be a $4$-dimensional manifold with a codimension-$2$ submanifold with isolated singularities.
We say that a $1$-form $\lambda$ is tame at  $U\subset \Sigma^{\reg}$ if on some \nbd $V$ of $U$  in $M$, there  
exists a smooth function $\kappa$ and a smooth $1$-form $\mu$, both with bounded derivatives, such that 
$$
\lambda=\kappa\alpha+\mu,
$$
where $\alpha$ is an angular form over $\Sigma^\reg$. 
\end{definition}
For simplicity, we will say that a $1$-form is tame at $\Sigma$ if it is tame at $\Sigma^{reg}$. 
The aim of this section is to produce 
tame Liouville forms associated to any polarization with reasonable singularities.
\begin{rk} Tameness is a differentiable notion : if the fibration or the connexion is modified, the class of tame forms remains unchanged.
\end{rk}
\begin{rk}\label{tameball} The Liouville forms defined in section \ref{liouvilleballsec} provides tame Liouville forms with arbitrary residues on the complement of a cross in a ball.
\end{rk}

\begin{rk} \label{tameradial}If $\lambda$ is a tame Liouville form with positive residues on $M\priv \Sigma$, where $\Sigma$ is a symplectic polarization,
the associated vector field $X_\lambda$ points outwards the regular part of $\Sigma$.
\end{rk}
\noindent {\it Proof :} Let $(z,R=r^2,\theta)$ be coordinates in a \nbd $V$ of a point $p\in \Sigma$, with 
$$
\begin{array}{l}
V\cap \Sigma=\{R=0\} \text{ and }
\om=dR\wedge d\theta+\pi^*\tau,
\end{array}
$$
where $\tau$ is a symplectic form on $\Sigma$ and $\pi(R,\theta,z)=(0,0,z)$. Tameness of $\lambda$ means that 
$$
\lambda=\kappa d\theta+\mu,\hspace{0,5cm} \mu\in \Gamma_V^\infty(\Lambda_1(T^*M)), \kappa\in \cc^\infty(V,\R),
$$
and positivity of the residue means that $\kappa$ takes positive values in $V$, provided $V$ is small enough. Thus,
$$
X_\lambda=\kappa\frac{\partial}{\partial R}+Z=\frac{\kappa}{2r}\frac{\partial}{\partial r}+Z,
$$
where $Z$ is a bounded vector field. Near $\Sigma=\{r=0\}$, this vector field is clearly radius increasing.\cqfd

\begin{rk}\label{rkglue} Tame forms are easy to glue. Namely, let $\lambda_1,\lambda_2$ be two tame $1$-forms, defined in 
\nbds $V_1,V_2$ of $U_1,U_2\subset \Sigma^{reg}$, with well-defined residues at $\Sigma$ (there might be several ones if $U_1$ or $U_2$ is not connected). Assume that the residues coincide, and that $\lambda_1-\lambda_2$ is exact on $(V_1\priv U_1)\cap (V_2\priv U_2)$. Then there exists a tame $1$-form $\lambda$ in $V_1\cup V_2$ which coincides with $\lambda_1$ 
on $V_1$ and with $\lambda_2$ on an arbitrary compact subset of $V_2\priv (\bar{V_1\cap V_2})$.
\end{rk}
\noindent {\it Proof :} By assumption, $\lambda_1-\lambda_2=dh$, with $h\in \cc^{\infty}(V_1\cap V_2\priv U_1\cap U_2)$. By the tameness assumption, 
$$
\lambda_1-\lambda_2(x)=[\kappa_1(x)-\kappa_2(x)]\alpha+\mu,
$$
where $\alpha$ is an angular form over $U_1\cap U_2$, $\kappa_1,\kappa_2$ are smooth functions and $\mu$ is a smooth
 $1$-form on
 $V_1\cap V_2$. Imposing the same residues implies that $\kappa_1-\kappa_2=0$ on $U_1\cap U_2$, so $\kappa_1-\kappa_2(x)=rb(x)$, where $r$ is a
 radial coordinate around $\Sigma$. Now the form $r\alpha$ is smooth on a \nbd of $\Sigma$, so $h$ is a smooth function on 
 $V_1\cap V_2$. This function can be smoothly extended by a function with arbitrary compact support in a \nbd of $V_1\cap V_2$. 
 The gluing is the form that coincides with $\lambda_1$ on $V_1$ and $\lambda_2+dh$ on $V_2$. \cqfd 
 
 \begin{lemma} \label{pdtame} If $\Sigma$ is smooth and $\sigma$ is a $2$-form representing the Poincaré-Dual class to $\Sigma$, 
 there exists a tame $1$-form $\lambda$ on $M\priv \Sigma$ with $d\lambda=\sigma$ and residue $1$ at $\Sigma$. The same holds when $\Sigma$ has isolated singularities, provided that each singularity has some \nbd already equipped with a tame $1$-form with residue $1$ at $\Sigma$.
 \end{lemma}
\noindent {\it Proof : } Consider first the case when $\Sigma$ is non-singular. Let $\cv$ be a tubular \nbd of $\Sigma$, $\pi:\cv\lra\Sigma$
a disc fibration and $\alpha$ an angular $1$-form associated to some connexion form on this fiber bundle. Inside $\cv$, $d\alpha$
represents  $\pd(\Sigma)$, so 
$$
d\alpha-\sigma=d\mu,
$$
where $\mu$ is a smooth $1$-form on $\cv$. The form $\alpha':=\alpha-\mu$ is therefore a tame primitive of $\sigma$, but is defined 
only in $\cv$. Choose now a $1$-form $\lambda$ on $M\priv \Sigma$ with $d\lambda=\sigma$. Let also $D$ be a local disc 
transverse to $\Sigma$ at $p$ and $\gamma_\eps$  a family of loops of $D$ converging to $p$. As in lemma \ref{extliouville} i), 
$\lambda(\gamma_\eps)$ has a limit, denoted by $\kappa$. Considering a deformation of $\Sigma$ whose intersections with 
$\Sigma$ are very close to $p$, and which coincides with a union of discs very close to $D$ near $p$, we see that $\kappa$ must be 
$1$, as soon as the normal bundle is non-trivial. On the other hand, when it is trivial, the $1$-form $\alpha'$ previously defined 
is closed. Hence, the $1$-form $\lambda-\kappa\alpha'$, defined on $\cv\priv \Sigma$, is closed. Thus,
$$
\lim_{\eps\to 0} \lambda-\kappa\alpha'(\gamma_\eps)=0,
$$ 
so $\lambda-\kappa\alpha'(\gamma_\eps)=0$ for all $\eps$. Now, if $\gamma_i$ is any basis of $H_1(\Sigma)$, there is a closed $1$-form
$\nu$ on $\Sigma$ with periods 
$$
\nu(\gamma_i)=\lambda-\kappa\alpha'(\wdt\gamma_i), 
$$
where $\wdt \gamma_i$ is a small perturbation of $\gamma_i$ in $\cv\priv \Sigma$. The form $\lambda-\kappa\alpha'-\pi^*\nu$
is therefore exact on $\cv\priv \Sigma$ so there is a function $h\in \cc^\infty(\cv\priv \Sigma)$ such that 
$$
\lambda-\kappa\alpha'-\pi^*\nu=dh.
$$
Any extension $\wdt h$ of $h$ to $M\priv \Sigma$ defines  a primitive of $\sigma$  on $M\priv \Sigma$ (by $\lambda-d\wdt h$), which is tame (because it coincides with $\kappa\alpha'-\pi^*\nu$ on $\cv$).

When $\Sigma$ has isolated singularities $(p_i)$, and when each singularity has a small \nbd $B_i$ equipped with a tame  $1$-form
$\lambda_i$ with residue $1$ at $\Sigma\cap B_i$, we proceed as follows. Notice first that the forms $\lambda_i$ can be assumed to 
be primitives of $\sigma$. Indeed, since $\lambda_i$ is tame and has a residue, the form $d\lambda_i$ extends to a smooth $2$-form
on $B_i\priv \{p_i\}$ (because $\lambda_i=\kappa_id\theta+rb(x)d\theta+\mu$ locally, and $rd\theta$ extends to a smooth $1$-form on $\Sigma^{reg}$). Since $H_2(B_i\priv \{p_i\})=0$, the form $\sigma-d\lambda_i=d\mu$ for a smooth one-form  $\mu$ defined on $B_i\priv \{p_i\}$. 
Correcting $\lambda_i$ by adding $\mu$, we therefore get a primitive of $\sigma$ on $B_i\priv \Sigma$, tame at $\Sigma\cap B_i$ and with residue $1$. Consider now a smooth perturbation $\Sigma'$ of $\Sigma$ which coincides with $\Sigma$ outside balls $B_i^\eps$
much smaller than $B_i$. By the previous argument, there exists a primitive $\lambda$ of $\sigma$ on $M\priv \Sigma'$, tame and with residue $1$ at $\Sigma'$. In view of remark \ref{rkglue}, we only need to understand that $\lambda-\lambda_i$ is exact in a 
\nbd of $\partial (B_i\priv \Sigma)$.  Since $\lambda$ and $\lambda_i$ have same residues, it amounts to showing that 
$\lambda-\lambda_i$ vanishes on the knot defined by $\partial B_i\cap \Sigma$. Consider a connected component $C$ of this knot, 
and fill it inside $B_i$ with a possible singular disc, which does not meet $B_i^\eps$, and whose singularities lie outside $\Sigma$. 
Then,
$$
\begin{array}{lll}
\lambda(C) & = \ds \int_D \sigma-D\cdot \Sigma' & \text{ (because $\lambda$ has residue $1$ at $\Sigma$)} \\
 & =\ds \int_D\sigma-D\cdot(\Sigma\cap B_i) & \text{ (because $D$ avoids $B_i^\eps$)}\\
  & =\lambda_i(C) & \text{ (because $\lambda_i$ has residue $1$ at $\Sigma$.}
\end{array}
$$
 \cqfd

\begin{cor}\label{tameliouville}
Let $(M,\om)$ be a symplectic manifold with a singular polarization $\Sigma=(\Sigma_i)$ with nodal singularities only. Then, if 
$$
[\om]=\sum \kappa_i\pd(\Sigma_i),
$$
there exists a tame Liouville form on $M\priv \Sigma$ with residues $a_i$ at $\Sigma_i$.
\end{cor}
\noindent {\it Proof :} Since each singularity consists of several branches intersecting transverally, there exists a tame 
Liouville form with residue 
$\kappa_i$ at $\Sigma_i$ near each singularity of $\Sigma_i$ (see remark \ref{tameball}). Consider now closed two-forms $\sigma_i$
on $M$ Poincaré-Duals to $\Sigma_i$, such that $\sum \kappa_i\sigma_i=\om$. By lemma \ref{pdtame}, there exists a tame $1$-form
$\lambda_i$ on $M\priv \Sigma_i$ with residue $1$ such that $d\lambda_i=\sigma_i$. Our tame Liouville form on $M\priv \Sigma$ 
is simply $\lambda:=\sum \kappa_i\lambda_i$.\cqfd  

Of course, in view of proposition \ref{liouvillell}, the same statement holds when the singularities have the form $\Pi_\alpha (z^p-\alpha w^q)$.
Let us now explain the relevance of tameness in our discussion. As we observed in remark \ref{tameradial}, the Liouville 
vector field associated to a tame form $\beta$ on $M\priv \Sigma$ is radial around the polarization - and it explodes in a controled way.
A basic consequence is that there is a well-defined $X_\beta$-trajectory emanating from any point of the polarization in each normal 
direction. In other terms, the flow of $X_\beta$ is well-defined on pairs $(p\in \Sigma,\theta)$ - where $\theta$ is a local angular 
coordinate around $\Sigma$ -, in other terms on the blow $\Sigma^*$ of $\Sigma$ (see \cite{moi4}, appendix A for details). As a result, we get the following (local) statement.  
\begin{prop}\label{locemb} Let $\lambda$ be a tame Liouville form on $B\priv \Delta :=B^4(1)\priv \{z_2=0\}\subset \C^2$, and 
$B_\lambda$ the basin of attraction of $\Delta$ in $B$ :
$$
B_\lambda:=\{p\in B\; |\; \exists \tau\, , \; \Phi_{X_\lambda}^{-\tau}(p)\in \Delta \}.
$$
Let also  $(M^4,\om,\Sigma,\beta)$ be a closed symplectic manifold with symplectic polarization $\Sigma$ and tame Liouville
form on $M\priv \Sigma$, and 
$$
\phi:\Delta^*\hra {\Sigma}^*
$$
be a lift to the blow-ups of an area preserving map between $\Delta$ and $\Sigma$.
Assume that the residues of $\lambda$ and $\beta$ at $\Delta$ and $\phi(\Delta)$ coincide. Then, 
 there exists a unique symplectic embedding 
$$
\Phi:B_\lambda(X)\hra M, 
$$
such that $\Phi_{|X}=\phi$ and $\Phi^*\beta=\lambda$.
\end{prop}

\section{Ball packings.} \label{ballsec}

\subsection{Proof of theorem \ref{packball}.}
\paragraph{Proof of i) :} Assume that $M$ has a symplectic packing by closed balls $B_i$ of capacities $a_i\in \Q$. Consider the 
blow-up $(\hat M,\hat \om)$ of $M$ along all the balls (recall that they are closed and disjoint). The cohomology class of $\hat \om$ is 
therefore 
$$
[\hat \om]=[\pi^*\om]-\sum_{i=1}^p a_i e_i,
$$
where $\pi:\hat M\to M$ is the blow-up map and the $e_i$ are the Poincaré-Duals of the exceptional divisors $E_i$ 
corresponding to the balls $B_i$. Since $\om$ is assumed to be a rational form, and the $a_i$ are rational numbers, the form 
$\hat \om$ is also rational. By Donaldson's result, for all $k$ sufficiently large and such that $k[\hat \om] \in H^2(M,\Z)$, there is a 
symplectic curve $\hat \Sigma_k$ Poincaré-Dual to $k[\hat \om]$. By definition, the homological intersection between $\Sigma_k$ 
and $E_i$ is $ka_i$, but there may be positive and negative intersections. The following lemma, which consists of an easy adaptation 
of Donaldson's proof  done in section \ref{donsec}, ensures that these intersections may be assumed to be positive. 
\begin{lemma}\label{donblowball}
For any $k$ sufficiently large and such that $k[\hat \om]\in H^2(M,\Z)$, there exists a symplectic curve $\hat \Sigma_k$ Poincaré-Dual 
to $k[\om]$ and whose intersections with the $E_i$ are transverse and positive. They therefore consist of $ka_i$ distinct points.
\end{lemma}
Now the pull-back of $\hat \Sigma_k$ by $\pi$ is a symplectic curve in $M\priv \cup B_i$ with boundaries. The boundaries are precisely
a union of $ka_i$ Hopf circles of each $\partial B_i$. Gluing to $\Sigma_k$ the corresponding Hopf discs inside the $B_i$, we get a 
symplectic curve which almost convenes. It is obviously  Poincaré-Dual to $k[\om]$, irreducible (because $\hat \Sigma_k$ is by Donaldson's method itself), and it has the desired singularities at the centers of the balls. The only problem is that $\Sigma_k$ may well be singular 
along the union of the Hopf circles if we do not care about the smoothness of the gluing (it is {\it a priori} only topological). However, in \cite{moi3}, lemma 5.2, it is proved that starting with the topological gluing, one can easily smoothen by a $\cc^0$-small perturbation 
localized inside the $B_i$ and near the Hopf circles. In short, the gluing can be smooth, so $\Sigma_k$ is the required curve. \cqfd

\noindent {\it Remark :} Lemma \ref{donblowball} can be easily proved using pseudo-holomorphic methods, as in \cite{mcpo}. Indeed, 
the curves $\Sigma_k$ given by Donaldson's method are $J_k$ holomorphic for some almost-complex structure $J_k$,  which are 
close to an almost-complex structure $J$ compatible with $\hat \om$. For instance, $J$ can be a structure inherited from a compatible 
almost complex structure on $M$  such that the $E_i$ are $J$-holomorphic (see \cite{mcsa}). Then, using automatic regularity of the exceptional 
$J$-spheres and their uniqueness (in a fixed homology class), we can deform the $E_i$ to exceptional $J_k$-spheres, which therefore 
intersect $\Sigma_k$ positively and transversally. By positivity of intersection, the $E_i'$ 
do not intersect, and these deformations come from a global Hamiltonian isotopy. Applying the inverse Hamiltonian flow to $\Sigma_k$, we 
get the needed symplectic curve, which proves lemma \ref{donblowball}. We present however a different proof in section \ref{donsec}, 
which only relies on Donaldson's techniques, and which is more suited to treat the singular situation raised by ellipsoid embeddings.

\paragraph{Proof of ii) :} Here is what we call the inflation procedure.
Consider a curve $\Sigma$ Poincaré-Dual to  $k\om$ with $p$ nodes, of mutliplicities $k_1a_1,\dots,k_pa_p$, with $k\leq k_i$ for all $i$. 
Assume that there are $p$ disjoint closed $k_ia_i$-crosses $X_i$ of size $a_i$, hence also disjoint crosses $X_i^\eps$ of size $a_i+\eps$. 
Each cross $X_i^\eps$ has a simply connected \nbd $\cv_i^\eps$ symplectomorphic 
(by a map $\phi_i$) to a 
\nbd $\cu_i^\eps$ of a union of $k_ia_i$ Hopf discs in $B(a_i+\eps)$. Consider the Liouville form $\lambda_i$ on $B(a_i)$ associated by proposition 
\ref{liouvilleball} to these $k_ia_i$ Hopf discs, with residues $\nicefrac{1}{k}$. It is also defined on $B(a_i+\eps)$, so since the closed crosses $X_i^\eps$  are disjoint, the form 
$$
\beta_{|\cv_i^\eps}:=\phi_{i*}\lambda_{i|\cu_i^\eps}
 $$
defines a Liouville form on $\cup \cv_i^\eps$ with residue $\nicefrac{1}{k}$. By lemma \ref{extliouville} iv), the restriction of $\beta$ to 
$\cv:=\cup \cv_i:=\cup (\cv_i^\eps\cap B(a_i))$ extends to a Liouville form on 
 $M\priv \Sigma$. Since $\nicefrac{k_ia_i}{k}\geq a_i$ for all $i$, proposition \ref{liouvilleball} iii) ensures that the (obviously disjoint) basins of attraction of the sets $\cv_i$, defined by 
 $$
 B_i:=\bigcup_{t\geq 0} \Phi^t_{X_\beta}(\cv_i),
 $$
  each  contain an open symplectic ball of capacity $a_i$. The embedding is given by 
 \fonction{\Phi_i}{B(a_i)}{B_i}{x}{\Phi_{X_\beta}^\tau\circ\phi_i\circ\Phi^{-\tau}_{X_{\lambda_i}}(x), \hspace{,5cm} \text{where }
 \Phi^{-\tau}_{X_{\lambda_i}}(x)\in \cu_i.} 
 Figure \ref{inflationfigure} is meant at illustrating this inflation process. 
 
\begin{figure}[h!]
\begin{center} 
\input inflation.pstex_t
\caption[Inflation]{Inflation of a \nbd of a cross.}
\label{inflationfigure}
{\footnotesize
\begin{itemize}
\item[(I)] $k_i< k$ :  inflating the cross does not lead to the embedding of the whole ball.
\item[(II)] $k_i=k$ :  inflating the cross allows the embedding of the whole {\it open} ball. 
                     Singularities on the boudaries are due to  the explosions of the Hopf circles
                           at the boundaries of the basin of attraction of the cross in the standard ball.
\item[(III)] $k_i>k$ :  inflating leads to an embedding of more than the ball. Then, the closed ball embeds smoothly, 
provided the crosses themselves are smooth.
\end{itemize} }
\end{center}
\end{figure}

 \subsection{Proof of theorem \ref{packopen}.}
Compared with the previous paragraph, we need to understand what happens when the  closures of the crosses are not disjoint 
any more. The argument is roughly the same, but is based on proposition \ref{locemb}, {\it i.e.} on  the analysis performed in the
 appendix to \cite{moi4}.  We first fix all the notations and describe the geometric picture.
In each $B(a_i)$, consider a cross $X'_i$ which  is symplectomorphic to $X_i$ ({\it i.e.} the k\"ahler angles between the discs are the same). Consider also disjoint symplectic balls $B_i^\eps$ in $M$, of capacity $\eps$, centered at the singular points $p_i$, and 
seen as the embeddings of the closed balls $B(\eps)\subset B(a_i)$ by a map $\phi_i$. We  further assume that these embeddings 
send $X_i'\cap B(\eps)$ to $X_i$. Denote by $\lambda_i$ the Liouville form on $B(a_i)$ with residues $\nicefrac{1}{k}$ given by proposition \ref{liouvilleball}. By corollary \ref{tameliouville}, we can extend the Liouville forms $\phi_*\lambda_i$ defined on  
$\cup B_i^\eps$ to a tame Liouville form $\beta$ on $M\priv \Sigma$. Our goal is now to extend the $\phi_i$ to embeddings of the 
open balls.\vspace*{,2cm}

\noindent {\bf Extension to ${X_i'}^*$.} Fix angular coordinates $\theta_i$ and $\theta'_i$ around $X_i\priv\{p_i\}$ and $X_i'\priv\{0\}$ respectively, and 
denote by $X_i^*, {X_i'}^*$ the blow-ups :
$$
X_i^*:=\{(p,\theta_i), p\in X_i\priv \{p_i\}, \theta_i\in \R\priv \Z\}.
$$
We can freely assume at this point that the restriction of $\phi_i$ to ${X_i'}^*\cap B(\eps)$ verify $\phi_i(p',\theta'_i)=(\phi_i(p'),\theta'_i)$. Consider now any extension $\wdt \phi_i$ of $\phi_{i|X'_i\cap B(\eps)}$ to an area-preserving map of $X'_i$ into $X_i$. 
Then, the formula 
$$
\phi_i(p',\theta'):=(\wdt \phi_i(p'),\theta_i) \hspace{,5cm} \text{ on } {X_i'}^*
$$
defines an extension of $\phi_i$ to ${X_i'}^*$ which is smooth on ${X_i'}^*$. 

 \begin{figure}[h!]
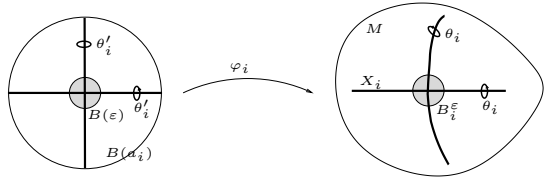

\begin{center} 
\input fatcross2.pstex_t
\caption{Blowing-up $\phi_i$.}
\label{fatcross}
\end{center}\vspace{-.7cm}
\end{figure}

\vspace{,2cm}

\noindent{\bf Inflation.} As we explained already, the tameness hypothesis implies that the positive flow of $X_\beta$ is well-defined on 
$X_i^*$. Similarly, the backward flow of a point $x\in B(a_i)$ either reaches $B(\eps)$ or a point on $X_i'$ (by proposition \ref{liouvilleball}, because $k_i\geq k$), along a certain direction, whose angle will be denoted $\theta'_i(x)$. We now inflate 
each $\phi_i$ to an embedding of the open ball $B(a_i)$ in the following way : 
\fonction{\Phi_i}{B(a_i)}{B_i}{x}{
 \left\{\begin{array}{ll} \Phi_{X_\beta}^\tau\circ\phi_i\circ\Phi^{-\tau}_{X_{\lambda_i}}(x), & \text{when }
 \Phi^{-\tau}_{X_{\lambda_i}}(x)\in  B(\eps),\\
 \Phi_{X_\beta}^\tau\circ\phi_i\big(\Phi^{-\tau}_{X_{\lambda_i}}(x),\theta_i'(x)\big), & \text{when }
 \Phi^{-\tau}_{X_{\lambda_i}}(x)\in  X_i'\priv B(\eps).
 \end{array}\right.
 } \vspace{,2cm}

The map is well-defined but this formula raises several questions. The first one concerns regularity : $\Phi_i$ is not 
obviously smooth, nor even continuous. However, $\Phi_i$ is clearly smooth on the basin of attraction of the open ball 
$B(\eps)$ and of the open annulus $X'_i\priv B(\eps)$ by proposition \ref{locemb}.  In order to see that $\Phi_i$ is actually 
smooth on $B(a_i)$, define $\Phi_i^2$ by the same formula, but replacing $B(\eps)$ by $B(\nicefrac{\eps}{2})$. The map 
$\phi_i^2$ is smooth on the basin of attraction of $X'_i\priv B(\nicefrac{\eps}{2})$ - which contains the locus where $\Phi_i$ 
is not known to be smooth. On the other hand, $\Phi_i^2$ coincide with $\Phi_i$ by the uniqueness part of proposition 
\ref{locemb}. Thus $\Phi_i$ is smooth on $B(a_i)$. The same proof shows that $\Phi_i$ is symplectic. Finally, the different 
embeddings are disjoint because they lie in the basin of attractions of disjoint subsets of the polarization.\cqfd

\subsection{Different shapes for efficient packings.}
The previous proof shows two important things. First, the natural objects that can be packed are the basins of attractions in $\C^2$ 
of crosses under the Liouville vector fields  associated to 
\begin{equation}\label{standgaliouville}
\left\{\begin{array}{l}
\lambda:=\frac{1}{n} \sum u_j^*\lambda_j, \text{ where }\\
\lambda_j:=(n\kappa_j-R_1)d\theta_1-R_2d\theta_2,\;   u_j\in U(2). 
\end{array}\right.
\end{equation}
 Moreover, the "vertical part of the embedding", in other terms the inflation process is smooth. If the discs of the crosses 
 are smooth up to the boundary, then the embedding of the corresponding basin of attraction is smooth up to the boundary, 
 locally near the boundary of the crosss. Combining these two observations, we get the following most general (and most ununderstandable) theorem. The notation $B\big( (\kappa_j)_{j=1\dots m}; (a_j)_{j=1\dots m}; (u_j)_{j=1\dots m}\big)$ stands 
 for the {\it open} basin of attraction in $\C^2$ of the "anisotropic" standard $m$-cross composed of standard open discs of sizes $(a_j)$ 
 inside complex lines $d_j:=u_j(\{z_2=0\})$ ($u_j\in U(2)$), under the flow of the Liouville vector field associated to  the form 
 given by (\ref{standgaliouville}). 
\begin{thm}\label{packgal}
Let ${\bf \Sigma}$ be a polarization of $M$
$$
[\om]=\sum_{l=1}^N \tau_l\pd(\Sigma_l).
$$
Assume that $\Sigma$ is covered by $r$ $m_i$-crosses each composed of $m_i^l$ discs of $\Sigma_l$, of sizes 
$(a_i^j)_{j\leq m_i}$, with respective positions described by $(u_i^j)_{j\leq m_i}$. Then $M$ has a full packing by 
$$
\coprod_{i=1}^r B\big((\underbrace{\tau_1,\dots,\tau_1}_{m_i^1},\dots,\underbrace{\tau_N,\dots,\tau_N}_{m_i^N}); (a_i^1,\dots,a_i^{m_i}); (u_i^j) \big).
$$ 
Moreover, around each smooth boundary point of the crosses in $\Sigma$, the corresponding embeddings are smooth up to the boundary. 
\end{thm}

This theorem may look a bit surprising, in that it seems to give too many different full packings. It is well-known and understandable
that curves of different degrees or different kind of singularities, which have no obvious symplectic relations, give different packings. 
But the parameters $(u_i)$ are too much. Indeed, easy local perturbations of the curves allow to modify freely those parameters. Why 
should it allow different embeddings ? The answer is that these basins of attraction do not really depend on these extra-parameters.
\begin{prop}\label{notwobasin} 
The basins $B\big( (\kappa_j); (a_j); (u_j)\big)$ and $B\big( (\kappa_j); (a_j); (v_j)\big)$ are domains of $\C^2$ with symplectomorphic boundaries.
\end{prop}
\noindent {\it Proof :} Consider the standard anisotropic crosses $X$ and $Y$ (associated to $\big( (\kappa_j); (a_j); (u_j)\big)$ and $\big( (\kappa_j); (a_j); (v_j)\big)$ respectively),
and a non-standard cross $Z\subset \C^2$ which coincides with $Y$ in $B(\eps)$ and with $X$ in the complement of $B(2\eps)$. 
Denote by $\wdt X$, $\wdt Y$ and $\wdt Z$ the infinite crosses in which $X$, $Y$ and $Z$ lie (just extend by complex lines). 
Let $\lambda_X$, $\lambda_Y$ be the standard Liouville forms on $\C^2\priv \wdt X$, $\C^2\priv \wdt Y$. Since $Y$ and $Z$ are 
symplectomorphic, there is a tame Liouville form $\lambda_Z$ on $\C^2\priv \wdt Z$ which coincides with $\lambda_Y$ inside $B(\eps)$ and with $\lambda_X$ outside $B(2\eps)$.  
Following now step by step the proof of 
theorem \ref{packopen}, we see that the basins of attraction of $X$ and $Y$ are symplectomorphic. On the other hand, the basin of attraction of $Z\priv B(2\eps)$ coincides with that of $X\priv B(2\eps)$.\cqfd

\subsection{Some explicit packings.}\label{maxpacksec}
We now explain how to construct geometrically maximal packings of $\P^2$ by $6$, $7$ and $8$ balls, proving theorem \ref{explicitballs}. 
\paragraph{Six balls :} Recall that the capacity of the optimal packing must be $\nicefrac{2}{5}$, and that the 
obstructions comes from $6$ conics passing each through the five of the  six centers of the balls. In fact, the 
same configuration of curves allows to construct the maximal packing. Indeed, fix six generic points $p_1,\dots,p_6$ in $\P^2$, and consider a conic through any five points of them.  The curve formed by the union of these conics has degree 
$12$, passes five times through each $p_i$, and it can be split into six crosses centered on the $p_i$, each disc of the cross having area $\nicefrac{2}{5}$ (simply divide each conic in five such discs). Moreover, the condition on the degree of this curve is verified : $\nicefrac{2}{5}\cdot 12\leq 5$. Thus, the six crosses can be inflated to a ball packing of maximal capacity. It may be worth noticing that  the produced packing is not very singular on its boundary : it can be choosed for instance to have exactly $10$ singularities at each balls, each singularity being rather nice (see the discussion in \cite{moi3} for a precise description of these singularities). They appear exactly on the polarization, when one parts the different rational curves into five discs. 

\paragraph{Seven balls :} In this case, the best capacity is $\nicefrac{3}{8}$ and the obstruction comes from a singular 
cubic passing once through six of the centers of the balls and twice through the last one. This particular cubic is not enough for our purpose of constructing seven balls of capacity $\nicefrac{3}{8}$, but almost. Consider seven generic points $(p_1,\dots,p_7)$ in $\P^2$. Consider one singular cubic through $(2p_1,p_2,\dots,p_7)$, one through 
$(p_1,2p_2,\dots,p_7)$ and one conic through $(p_3,\dots,p_7)$. The union of these curves is a curve of degree $8$ that passes exactly $3$ times through each $p_i$. Split each cubic into $8$ discs of area $\nicefrac{3}{8}$ centered
on the $p_i$ and consider also discs of area $\nicefrac{3}{8}$ around $p_3,\dots,p_7$ in the conic (it is possible because 
$5\cdot\nicefrac{3}{8}\leq 2$). We thus have found a cross of area $\nicefrac{3}{8}$ with $3$ branches around each $p_i$ inscribed in a curve of degree $8$. That is what is needed for getting the optimal packing. Regularity on the boundary is however more difficult to handle, precisely because in this situation, the capacity of the balls ($\nicefrac{3}{8}$) is really the same as the quantity multiplicity/degree. Thus the singularities that appear in the partition of the curves into discs propagates along the Liouville vector field (hence giving rise to "nice" discs of singularities), and new ones appear (which we cannot control at all), at the Hopf discs that are not in the basin of attraction of the cross. 

\paragraph{Eight balls :} The best capacity is $\nicefrac{6}{17}$. In order to produce the packing, we therefore need a curve of degree $17$ passing six times through each points, or a curve of degree $34$ passing $12$ times through each points and so on. In fact we will produce a curve of degree $51$ passing $18$ times through each point. The obstruction comes from a sextic, passing three times through one point, and twice through all other points. As before, the curve that 
allows the construction will be obtained from this "obstruction" curve. Namely, the curve we consider is composed 
of eight sextics, each passing twice through each point but one (each time different), where it passes three times. This curve has degree $48$ and passes exactly $17$ times through each point. Add to this curve a cubic passing once through each of the eight points, and you get the announced curve of degree $51$. 

\subsection{Proof of theorem \ref{packballirr}.}
First assume that there are $p$ disjoint closed balls $B_i$ of size $(a_i)_{i=1\dots p}$ in $M$. Consider rational symplectic 
forms $\om_1,\dots, \om_N$ very close to $\om$, such that each $\om_l$ coincides with $\om$ in a \nbd of the balls, and such that 
$[\om]$ lies in the convex hull of the $[\om_l]$ in $H^2(M,\R)$ :
$$
[\om]=\sum_{l=1}^N \mu_l [\om_l], \hspace{,5cm} \sum \mu_l=1,\; \mu_l>0.
$$ 
Since the $\om_l$ coincide in a \nbd of the balls, the symplectic blow-up of $(M,\om_l)$ along the $p$ balls provide different 
symplectic forms $\hat \om_l$ on the same manifold $\hat M$, with the same exceptional divisors $E_i$. Now, applying lemma
\ref{donblowball}, we find $\hat \om_l$-symplectic curves $\hat \Sigma_l$ Poincaré-Duals to $k_l\hat \om_l$ and whose intersections 
with the $E_i$ are transverse and positive. It will also follow from the proof of this lemma that the intersections of $E_i$ 
with the different $\Sigma_l$ can be assumed to be disjoint (see remark \ref{disjointintersection}). In fact, we can even require 
that the $\hat \Sigma_l$ are $\om$-symplectic, and intersect positively and transversally in $\hat M$ (see \cite{moi5}, theorem 2). 
Projecting these curves down to $M$, and gluing smoothly the convenient Hopf discs, we get curves $\Sigma_l$ Poincaré-Duals 
to $k_l\om_l$, and whose intersections with $B_i$ consists of $k_l^i=k_la_i$ Hopf discs. Now,
 ${\bf \Sigma}:=(\Sigma_l,\nicefrac{\mu_l}{k_l})$ is a polarization of $M$ with $p$ nodes of multiplicities $k^i:=\sum k_l^i$, 
 which contains the $p$ disjoint $k^i$-crosses of size $a_i$ by construction. Finally, 
 $$
 \sum_l k^i_l\tau_l =\sum_l k_l^i\frac{\mu_l}{k_l}=\sum_l k_la_i\frac{\mu_l}{k_l}=\sum_l \mu_l a_i=a_i.
 $$
 If you wish to verify that the left hand side above may even be strictly larger than $a_i$, simply notice that $M$ contains 
 also disjoint closed balls of radii $a_i+\eps$. 
 
 Conversely, the proof that a polarization with the properties listed in theorem \ref{packballirr} allows to construct 
the ball packing  is exactly  the same as in theorem \ref{packball} ii). \cqfd

\section{Ellipsoid embeddings.} \label{ellsec}
As for balls, it will be convenient to be able to blow-up the ellipsoids. The result is a symplectic orbifold, that we describe now. 
The reader can also consult \cite{godinho,nipa}.

\subsection{Blowing-up an ellipsoid.}
\paragraph{Weighted projective space.} In this paragraph, we fix two relatively prime integers $p,q$. It is well-known 
that $\P^2$ is symplectically a closed ball, whose boundary has undergone a symplectic reduction, {\it i.e.} has been collapsed along its characteristic foliation. We define similarly $\P^2(1,p,q)$ as the ellipsoid $E(p,q)$ whose boundary 
has been collapsed along its characteristic foliation. It is easy to see that this definition gives the same topology to 
$\P^2(1,p,q)$ as the more classical one 
$$
\P^2(1,p,q):=\C^3_{/\C}, \text{ with action } \xi\cdot z:=(\xi z_1,\xi^p z_2,\xi^q z_3), 
$$
but it has the advantage of giving an easier geometric model for the symplectic structure. Notice that it is easy to see 
(using any of these two models) that the divisor at infinity $\{z_1=0\}$, which corresponds to the boundary of the ellipsoid after the collapsing process, is topologically the weighted projective space $\P^1(p,q)$ - we determine its symplectic area below.

Let $\Phi:B(pq)\lra E(p,q)$ be the ramified $pq$-covering of the ellipsoid by the ball described in the previous paragraph. It obviously sends the characteristic foliation of $\partial B$ to the one of $\partial E(p,q)$, so it descends to a symplectic 
(ramified) covering of $\P^2(1,p,q)$ by $\P^2(pq)$ (meaning that the projective line has area $pq$). The importance of this covering is two-fold. First, it explicitly presents $\P^2(1,p,q)$ as a symplectic orbifold, with an explicit desingularizing 
map. Notice that the divisor at infinity, which is itself singular, is also desingularized by this map, since it corresponds to a projective line. The second point is that it allows to complete the description of our weighted blow-up by giving us the symplectic area of the divisor at infinity. It is precisely $1$, since it is a quotient of a symplectic line of area $pq$ by a symplectic group covering of degree $pq$. Let us sum up this discussion :
\begin{prop}\label{wbu}
Let $p,q$ be two relatively prime integers. The weighted projective space $\P^2(1,p,q)$ is obtained from the 
ellipsoid $E(p,q)$ by collapsing the Hopf fibration of $\partial E(p,q)$. It is a symplectic orbifold, with group $\Z_p\times \Z_q$, and a desingularization map $\Phi:\P^2(pq)\lra\P^2(1,p,q)$. Finally, it has three "distinguished" curves : the horizontal one of area $p$, the vertical one of area $q$, and the divisor at infinity, of area $1$.
\end{prop}
 
 \paragraph{Blowing-up ellipsoids.} Assume now that a closed ellipsoid $E(a,b)=\tau E(p,q)$ embeds into a symplectic manifold $(M,\om)$. 
 \begin{prop}\label{orbiblow}
 Removing from $M$ the interior of the ellipsoid $E(a,b)$ and collapsing the characteristic foliation of the resulting boundary (which is $\partial E(a,b)$), we get a symplectic orbifold with one or two singularities located on the exceptional divisor (the projection of $\partial E(a,b)$). This 
 exceptional divisor is a symplectic suborbifold of area $\tau$, which is desingularized by a symplectic desingularization of  $M$.
 \end{prop} 
The resulting symplectic orbifold will be called the blow-up of $M$ along $E(a,b)$ and denoted $(\hat M,\hat \om)$. Recall that classical 
blow-up gives a presentation of a symplectic manifold with some ball inside as a Gompf sum of the classical blow-up and the projective 
space along the exceptional divisor on one side and a projective line on the other. Similarly, this singular notion of blow-up allows to 
think at a symplectic manifold with some ellipsoid $E(a,b)$ inside as a Gompf sum of two symplectic orbifolds - the blow-up and the weighted projective space $\tau \P^2(1,p,q)$ - along symplectic suborbifolds - the exceptional divisor and the line at infinity. 

\noindent {\it Proof of proposition \ref{orbiblow} :} Observe that the map $\Phi:\tau B(pq)\to E(a,b)$ defined above extends to the whole 
of $\C^2$, hence in a \nbd of $\tau B(pq)$. Moreover, if the closed ellipsoid $\bar E(a,b)$ symplectically embeds into $M$, so does a slightly 
larger ellipsoid. The extension of $\Phi$ to a covering of this larger ellipsoid therefore gives a symplectic desingularization of the manifold 
described in proposition \ref{orbiblow}.\cqfd

\subsection{Proof of theorem \ref{packell}.}
Let us first assume that $\ce:=E(a,b)=\tau E(p,q)$  embeds into $(M,\om)$ where $\om$ is a rational class ($p,q$ are relatively prime integers). Consider the blow-up $\hat M$ of $M$ along this ellipsoid, and call $E$ the resulting exceptional divisor. Provided that $\tau$ is rational, the symplectic form $\hat \om$ on $\hat M$ is rational. Although the manifold is singular and the present framework is that of symplectic 
orbifolds, the next lemma asserts that Donaldson's techniques can be carried out in this setting.
\begin{lemma} \label{donblowell}
For any $k$ sufficiently large and such that $k[\hat \om]\in H^2(\hat M,\Z)$, there exists  symplectic smooth curves $\hat \Sigma_k$ 
 Poincaré-Duals to $k\hat \om$, which intersects the exceptional divisor transversally and positively at 
 exactly $k\tau$ regular points.
\end{lemma}
\noindent Projecting a curve $\hat \Sigma_k$ down to $M$, we get  a symplectic curve $\Sigma_k$
 whose boundary is made of $k\tau$ regular characteristic leaves of $\partial \ce$.  Now each such characteristic 
 bounds a complex curve of equation $z_1^q=\alpha z_2^p$ in $\ce$. Gluing these complex curves to $\Sigma_k$, we get a 
 closed curve Poincaré-Dual to $k\om$. As before, this gluing can be assumed to be smooth by \cite{moi3}. The resulting symplectic 
 curves has therefore all the required properties. 
 
 \paragraph{} Conversely, assume that $M$ contains a symplectic curve $\Sigma$ Poincaré-Dual to $k\om$, with a 
 singularity composed of $k\tau$ branches locally of the form $z_1^q=\alpha z_2^p$. Namely, we assume that 
 a \nbd of the singularity has a Darboux chart where $\Sigma$ restricts to the curve 
 $$
 \prod_{i=1}^{k\tau} (z_1^q-\alpha_i z_2^p)=0.
 $$
Now, cut the part of this curve that lies inside a very small ellipsoid $\tau' E(p,q)$ centered at the singularity $x$, 
and glue back smoothly the cone over the resulting boundary, with center at the  origin.  The curve $\Sigma'$ that results 
from this operation is Poincaré-Dual to $k\om$ and locally modelled to a cone singularity over a $(p,q)$-torus link. Assuming 
that the volume obstruction is satisfied, we get :
$$
k\tau\cdot \tau pq=kab\leq k\vol M=\int_M k\om\wedge\om \overset{\pd}{=} \int_\Sigma \om=\ca_\om(\Sigma),
$$
so $\Sigma$ contains $k\tau$ disjoint (singular) discs of area $\tau pq$ that meet at $x$, that is a 
 copy of  $\Upsilon:=\cup_{i=1,\dots,k\tau}\Delta_{\alpha_i}(p,q)$ (see definition \ref{defcone}). Now, 
$M\priv \Sigma$ has  a Liouville form with residue $\nicefrac{1}{k}$ at $\Sigma$ as well as $E(a,b)\priv \Upsilon$. 
The associated Liouville vector field on $E(a,b)$ is radius increasing by \ref{liouvillell} ii) because 
$k\tau\cdot \nicefrac{1}{k}\geq \tau$, so the same construction as for theorem \ref{packball} embeds the 
ellipsoid $E(a,b)$.\cqfd

\section{Donaldson curves in blow-ups.} \label{donsec}
 \subsection{In the blow-up of a ball.}\label{donballpart2}
We prove now lemma \ref{donblowball} avoiding the use of pseudo-holomorphic curves. The general statement is the following :
\begin{prop}\label{donsub}
Let $(M,\om)$ be a symplectic manifold with $\om\in H^2(M,\Z)$. Let $N$ be a closed symplectic submanifold of 
arbitrary codimension. Then, for all sufficiently large $k$, there exists a symplectic hypersurface $\Sigma$ 
Poincaré-dual to $k\om$ and such that $\Sigma\cap N$ is a symplectic hypersurface of $N$ 
Poincaré-Dual to $k\om_{|N}$.
\end{prop}
Applying this proposition to the blow-up and its exceptional divisor yields a proof of lemma \ref{donblowball}. 
Although it is easy and rather folkloric (a more subtle version for Lagrangian submanifold can be found in \cite{augamo}), we 
sketch a proof of proposition \ref{donsub}, because we also need to adapt it in the less usual setting of symplectic orbifolds.

\paragraph{Quick review of Donaldson's construction.} Let $J$ be a compatible almost complex structure, $\cl$ a hermitian line bundle over $M$ with curvature $\om$. Donaldson's proof of the existence of a symplectic hypersurface consists in producing an approximately holomorphic and $\eta$-transverse section of $\cl^{\otimes k}$. This section 
is obtained as a perturbation of the (obviously holomorphic) zero-section by approximately-holomorphic peak-sections 
very localized around the points of the manifolds. These peak-sections decrease exponentially fast around a point, 
and have support in a ball of radius of order $\nicefrac{1}{\sqrt[3] k}$. Let us describe a bit more the peak-section $\sigma_p^k$ around a point $p\in M$. Identify a \nbd of $p$ in $M$ to a ball in $\C^n$ by a map $\chi_p$ that takes 
$p$ to the origin, $\om$ to the standard symplectic form on $\C^n$ and $J$ to $i+O(|z|)$. Inside this ball, the 
bundle $\cl^{\otimes k}$ has a connection given by $d+k\sum z_jd\bar{z_j}-\bar{z_j}dz_j$ and in the trivialization of the 
bundle given by this connection, the peak section is 
$\sigma_p^k=\chi(z)e^{-k|z|^2}$, where $\chi(z)$ is a cut-off function at the right scale (notice that $e^{-k|z|^2}$ is $i$-holomorphic). Donaldson method consists 
in adding inductively to the zero-section some small enough multiple of these peak-functions, centered in points on a regular grid, so that the perturbation gain transversality on the balls where the peak-section is added, while being cautious enough not to destroy the transversality already obtained.

\paragraph{Proof of lemma \ref{donsub}.} Take $J$ compatible with $\om$, but also to $\om_{|N}$ (in particular, $J$
preserves $TN$). Notice that the restriction of the bundle $\cl^{\otimes k}$ to $N$ has curvature $k\om_{|N}$. 
By classical \nbd theorems, a \nbd of a point $p\in N$ can be identified symplectically to a symplectic 
ball, where  $J$ is again $i+O(z)$ and $N$ is a (complex) linear subspace. In these coordinates, the above mentioned peak-section
$\sigma_p^k$ on $M$ restricts to the peak section $\sigma_{N,p}^k$ on $N$ (associated to $(N,\om_{|N},J_{|N})$). Therefore, perturbing the zero section 
first above $N$ by Donaldson's receipe provides a section which is zero at distance $\nicefrac{1}{\sqrt[3] k}$ of 
$N$, and which is $\eta$-transverse to zero over $N$. If we then continue the process with perturbation sufficiently 
small, we do not affect the transversality of the section over $N$, which implies that the intersection of the vanishing locus of the final section with $N$ is precisely a symplectic polarization of degree $k$ of $N$. \cqfd

For the proof of theorem \ref{packballirr}, we need a slightly stronger version of this statement. Namely, we must know that if 
$\om_1,\dots, \om_m$ are rational and close enough to a given symplectic form $\om$, the intersections of Donaldson's hypersurfaces 
$\Sigma_1^k,\dots,\Sigma_m^k$ with $N$, which are Donaldson's hypersurfaces of $N$ as we already saw, can be assumed to intersect 
transversally (in dimension $4$, it will imply that these intersections are disjoint). This point is obtained by requiring the $\eta$-transversality 
for the sections of $\cl_i^{\otimes k}\oplus \cl_i^{\otimes k}$ (where $\cl_i$ is Donaldson's bundle associated to $\om_i$). This analysis 
is precisely done in \cite{moi5}. As in the previous proof, getting this transversality first above $N$ and then on $M$ gives the desired result.
\begin{rk}\label{disjointintersection}
In dimension $4$, the intersections of $\om_l$-Donaldson's hypersurfaces with a symplectic curve $N$ can be assumed to 
be disjoint, provided the $\om_l$ are close enough to a fix symplectic form.
\end{rk}

\subsection{In the blow-up of an ellispoid.}
In this paragraph, we explain that Donaldson's techniques  generalize to the setting of $4$-dimensional orbifolds, where lemma \ref{donsub} has the following analogue :
\begin{lemma}\label{donsuborbi}
Let $(M,\om)$ be a 4-dimensional symplectic orbifold with isolated cyclic singularities and $\om\in H^2(M,\Z)$. Let $N$ be a closed symplectic suborbifold which is desingularized in the uniformizing charts of $M$ (in particular, $N$ has only isolated cyclic singularities). Then, for all sufficiently large $k$, there exists a symplectic curve $\Sigma$ Poincaré-dual to $k\om$, which avoids all singular points of $M$ and $N$, and whose intersection with $N$ is transverse and positive. 
\end{lemma}
 This lemma may well be true in a more general setting but this version is enough for proving lemma \ref{donblowell}, which is our prupose. 
 Before proving this lemma, we first need to understand that Donaldson's method works in symplectic orbifolds.
 
 \begin{thm}\label{donorbi} Let $(M,\om)$ be a $4$-dimensional symplectic orbifold with isolated cyclic singularities, and $[\om]\in H^2(M,\Z)$. 
 Then for all sufficiently large $k$, there is a smooth symplectic curve $\Sigma_k$ Poincaré-Dual to $k[\om]$ which avoids all singular 
 points of $M$. 
 \end{thm}
\noindent {\it Proof :} As explained in the previous paragraph, Donaldson's method is based on several ingredients. One needs :
\begin{itemize}
\item[($1$)] A compatible almost complex structure $J$ on $(M,\om)$ (call $g$ the induced metric), 
\item[($2$)] Regular grids $\Gamma_k$ of $M$ at scale $\nicefrac{1}{\sqrt k}$,
\item[($3$)] A line bundle $\cl$ on $M$ of curvature $\om$,
\item[($4$)] The highly localized approximately holomorphic sections $\sigma_p^k$.
\end{itemize}
Then the method consist in adding iteratively to the zero-section (or any approximately holomorphic sections) of $\cl^{\otimes k}$
some linear combination of the $\sigma_p^k$, where $p$ belongs to the grid $\Gamma_k$. The first question is to decide whether the 
ingredients of the receipe are available in the orbiworld. An almost complex structure is easy to produce : in a uniformizing chart near the singularity, choose $J=i$. Since the action of the covering group is by unitary maps, it preserves $i$ so $J$ defines a complex structure on 
some \nbd of the singularities, which we can extend to an almost complex structure by general arguments. Similarly, in the uniformizing chart, 
the connexion $d+\sum z_jd\bar z_j - \bar z_jdz_j$ (with curvature $\om$) is also unitary-invariant, so it projects to a connexion 
on the trivial line bundle over a \nbd of the singularities with the good curvature. Extending the fiber bundle to the whole of $M$ then
goes as in the smooth case. The regular grid is also not a problem. But the fourth object above mentioned has no obvious equivalent. 
Indeed,  the controled decay rate for the peak sections depends heavily on the existence of Darboux charts around each point of the 
manifold which depend smoothly on the parameter. But this point fails dramatically on symplectic orbifolds, because Darboux charts 
around regular points close to a singularity cannot contain this singularity, so they must be very distorted. Hence, this straightforward 
approach does not go through. Instead we proceed as follows. We first fix a Darboux uniformizing chart around each singularity. Namely, 
this is a symplectic covering 
$$
\Phi:B^4(1)\lra M
$$
ramified along the axis $\{z_i=0\}$, and invariant under the symplectic action of $G:=\Z_p\times \Z_q$ on $B^4(1)$ given 
by 
\begin{equation}\label{eqaction}
(\xi_1,\xi_2)\cdot (z_1,z_2)=(\xi_1z_1,\xi_2z_2)
\end{equation}
These uniformizing charts also come with a complex structure $i$ and a line bundle $\wdt \cl$ with curvature $\om$, 
which projects to $\cl$. The main point is that we can produce good sections of $\wdt \cl^{\otimes k}$ over $B(1)$ which are
 invariant by the covering group $G$. 
\begin{lemma}\label{donballgroup}
There exists a sequence $(s_k)$ of sections of $\wdt \cl^{\otimes k}$ over $B(1)$ which are compactly supported, 
approximately $i$-holomorphic, $\eta$-transverse on $B(1-\eps)$ (for fixed $\eps,\eta$), non-vanishing at the origin, 
and invariant by the action of $G:=\Z_p\times \Z_q$ given by formula (\ref{eqaction}). 
\end{lemma}
Unfortunately, this lemma does not seem to follow formally from Donaldson's construction, so we explain the relevant modifications 
below. Notice before that it obviously concludes the proof, because these equivariant sections on $B(1)$ define sections on fixed \nbds of the 
singularities that can be extended (by $0$) to sections $s_k$ of $\cl^{\otimes k}$ which are approximately holomorphic on $M$,
$\eta$-transverse on some fixed \nbd of the singularities, and which do not vanish at the singularities. Applying Donaldson's 
perturbative method to $s_k$ outside this \nbd, we get a section of $\cl^{\otimes k}$ whose zero-set has the required properties,
provided $k$ is large enough.\cqfd
\noindent {\it Proof of lemma \ref{donsuborbi} :} It will be clear from the proof of lemma \ref{donballgroup} that when we
are given an additional data of a $G$-invariant set of symplectic submanifolds in $B^4(1)$, we can also require that 
the $G$-invariant sections given by lemma \ref{donballgroup} are approximately holomorphic and $\eta$-transverse when 
restricted to these submanifolds. The induced section (on the \nbd of the singularities) can be extended by $0$ on $M$, and the
procedure explained in section \ref{donballpart2} can be applied in the rest of the manifold (which is non-singular) in order 
to produce the section whose zero-locus is the desired symplectic curve. \cqfd

\subsection{Proof of lemma \ref{donballgroup}.}
Before getting to the actual proof, let us observe that the naive approach of averaging a (non-equivariant) Donaldson's 
section under the action of $G$ may not produce the desired result, since interference between the different terms of the 
sum may destroy the transversality. The strategy of the proof is therefore to follow closely Donaldson's argument, and make
it equivariant step by step. 

Let us first recall how Donaldson's method produces a sequence of sections which may check every condition of lemma \ref{donballgroup},
except equivariance. To the zero section of the line bundles $\cl^{\otimes k}$, we add inductively the localized peak sections $\sigma_p^k$
introduced in section \ref{donballpart2}, centered at points of the grid $\Gamma_k$. For instance, in our standard situation of euclidean 
metric,  the grid can be taken to be a regular lattice of girth $\nicefrac{1}{2\sqrt k}$ (so that the balls of radii $\nicefrac{1}{\sqrt k}$ cover 
$B(1)$). This lattice can be assumed to be invariant by the group action, as well as the different peak functions associated to these points. 
Now, the  iteration goes as follows. The  grid $\Gamma_k$ is parted into subgrids $\Gamma_k^1,\dots, \Gamma_k^N$ (we sometimes refer to $[1,N]$ as to the colors of the points). We first add to the zero-section a linear combination 
$$
\sum_{p\in \Gamma_k^1} w_1(p)\sigma_{p}^k,
$$
where  the $w_1(p_i)\in \C$ have magnitude $1$, and can be choosen independently. In particular, if $\Gamma_k^1$
is invariant by the group action, we can freely decide that $w_1$ take the same value for points in the same orbit (because we perturb the 
zero-section which is equivariant). The new section $s_1^k$ is therefore equivariant and $1$-transverse on a \nbd of size 
$\nicefrac{1}{\sqrt k}$  of $\Gamma_k^1$. Now the next step consists in adding to $s_1^k$ a linear combination 
$$
\sum_{p\in \Gamma_k^2} w_2(p) \Sigma_{p}^k,
$$
where $w_2(p)$ are this time of order much smaller than $\min{w_1(p)}\approx 1$.
Again, provided that $\Gamma_k^2$ is $G$-invariant,  this perturbation can lead to an equivariant section $\eta(\alpha)$-transverse. We 
then  repeat the process, and the analysis performed by Donaldson shows that the number $N$ of subgrids (equivalently the number of 
steps of the induction) can be taken not too large, so that the method ends-up with $\eta_N$-transversality, where $\eta_N$ is bounded 
from below (with respect to $k$). 

We thus conclude that if the coloring can be made $G$-invariant, lemma \ref{donballgroup} is proved by Donaldson's line of arguments.
Now the condition on the coloring is that any two points of the same color must be sufficiently far away so that the peak sections 
centered at these two points do not interfere and destroy one another  the transversality they achieve individually. Since 
the peak-sections have norm 
$$
|\sigma_p^k(q)|\approx e^{-kd(p,q)^2},
$$
we see that the suitable condition is that two points can have the same color as soon as they are at least $\nicefrac{D}{\sqrt k}$-distant, 
where $D$ is some large enough constant. Since the orbits of a point far from the origin are well-separated, we get :
\begin{cor}
The set $\Gamma_k^i\priv B(0,\nicefrac{D'}{\sqrt k})$ (where $D'=|G|D$) can be choosen to be $G$-invariant.
\end{cor}
Thus, the coloring can be made equivariant if and only if the origin is the only point of $\Gamma_k$ in $B(0,\nicefrac{D'}{\sqrt k})$. 
We finally need to understand that the grid can be distorted near the origin  so that it checks this condition (recall we initially asked for 
a $\nicefrac{1}{2\sqrt k}$ net). Consider now the section
$$
s_0^k:=e^{D'}\sigma_0^k=e^{D'}e^{-k|z|^2}\chi(z), 
$$
where $\chi(z)$ is a cut-off function which is $1$ until $B(1-\eps)$. This section is {\it holomorphic} and $1$-transverse on 
$B(0,\nicefrac{D'}{\sqrt k})$, and approximately holomorphic on $B(1)$. At this point, it is very important to choose $J=i$ at least on a 
\nbd of the  origin, 
because if $\sigma_0^k$ had only been approximately holomorphic, the factor $e^{D'}$ would have altered the estimate and the 
argument may well have failed. To conclude, instead of starting Donaldson's process with the zero section, simply start with the 
section $s_0^k$, which only needs corrections outside $B(\nicefrac{D'}{\sqrt k})$ to achieve transversality. \cqfd

\section{Miscallenous remarks.}
\paragraph{Gromov's non-squeezing.} Since the criterion of theorem \ref{packball} for the embedding problem only involves 
the existence of  (almost) Donaldson curves with some prescribed singularities, it is natural to ask whether Gromov's 
non-squeezing theorem can follow from Donaldson's techniques instead of pseudo-holomorphic curves 
(at least in dimension 4), which would be highly unexpected. Namely the question is the following :
\begin{question} Let $\Sigma$ be a symplectic curve of $(S^2\times S^2,\om:=\om_0\oplus2\om_0)$ which is Poincaré-Dual to 
$k\om$. If $\Sigma$ has  exactly one singularity, on a point $p$, locally modeled on $N$ holomorphic discs intersecting transversally, 
can one prove that $N\leq k$.
\end{question} 
The answer is obviously yes in the framework of $J$-holomorphic curves. Indeed, the curve $\Sigma$ is $J$-holomorphic 
for some almost complex structure, and by Gromov's argument, we can find a  $J$-curve $S$ in the class of $S^2\times \{*\}$ passing through $p$. Then, by positivity of intersections :
$$
1\cdot k+0\cdot 2k=k=\Sigma\cdot S\geq N\cdot 1 \hspace{0,5cm}\text{(the local intersection).}
$$
Although Donaldson's method sometimes allows to deal with positivity of intersections, it seems hardly useful 
 to produce a curve in this class of homology (or even a multiple of this class, or something close).  

\paragraph{Minimal degree of curves with prescribed singularities \cite{grlosh,shustin}.} 
This paper is about relations between singular curves and symplectic packings in dimension $4$. In section \ref{maxpacksec}, we use 
 algebraic geometry to produce symplectic packings. But the reverse direction can also be investigated.
\begin{question}
Given a list of singularity types $(S_i)_{i=1\dots r}$, what is the minimal degree $d(S_i)$ of an algebraic curve of $\P^2$ with 
$r$ singular points of type $S_i$ ?
\end{question}
Of course, "singularity type" can have various meanings (up to analytic, or topological equivalence), which may lead to different answers. 
We refer for instance to \cite{grlosh,shustin} for algebraic approaches to this problem. Now the remark is the following : 
in view of theorems \ref{packball1} and \ref{packell}, and since the packing problem was completely solved for any configuration 
of balls and ellipsoids in $\P^2$ by McDuff in \cite{mcduff4}, the corresponding question in a symplectic setting 
 has a complete solution, at least asymptotically in the degree. If one believes in the symplectic isotopy conjecture, it gives the optimal result. In any case, any symplectic lower bound for $d(S_i)$ is also an algebraic one, and it is sometimes better than what is
 known. For instance, we know from \cite{mcduff4} that the closed ellipsoid $\tau E(1,6)$ embeds into $\P^2$ if and only if 
 $\tau\leq \nicefrac{2}{5}$. Thus, 
 \begin{thm}
 Let $S$ be a singularity with a symplectic model given by $m$ smooth branches intersecting in one point with tangency order $6$ :
 $$
 S:=\left[\{w(w-z^6)\cdots(w-(m-1)z^6)=0\}\right].
 $$ 
 Then, $d(S)\geq \frac{5}{2}m$. 
 \end{thm} 
Notice that the classical bound $\sqrt {\mu(S)}+1=\sqrt{6m^2-7m+1}+1\approx 2,45m$ is smaller that $2,5m$ when $m$ is large.

\paragraph{Isotopies of balls.} As in \cite{moi4}, it may be tempting to use theorem \ref{packball} to try to isotop balls in a
symplectic $4$-manifold. The idea is that an embedding of a ball is determined by a symplectic curve with a prescribed 
type of singularities, together with a Liouville form on the complement of this symplectic curve. Now, since Liouville forms are 
soft objects (they form a convex set for instance), one could easily prove that some balls are isotopic, provided he knows such 
an isotopy result on symplectic curves with some prescribed singularities. In particular, if the problem of finding a curve in the 
homology class of $\pd(k[\om])$ with some nodal singularity has a non-vanishing Gromov-Witten invariant, it seems likely that 
one can prove an isotopy results for the balls of the maximal size that can be produced from this curve. This approach is slightly 
more direct, but is essentially the same as McDuff's \cite{mcduff3,mcduff4}. 

\paragraph{McDuff's result on ellipsoid embeddings.}
Let us recall that McDuff proved in \cite{mcduff4} that on some manifold, the embedding problem for an ellipsoid with fixed shape 
is equivalent to the packing problem for some balls with given ratio between the different radii. As in the previous paragraph, it seems 
natural to ask whether this remark also proceeds from the results contained in this paper. In some sense, the question is whether 
when we have a Gromov-Witten invariant for some problem involving a curve in some class with some simple singularities (several 
nodal points for instance), one can deform this curve in order to produce a new curve with one more complicated singularity 
(such as considered in theorem \ref{packell}).

{\footnotesize
\bibliographystyle{abbrv}
\bibliography{bib3.bib}
}

\vspace{2cm}
\noindent Emmanuel Opshtein,\\
Institut de Recherche Mathématique Avancée \\
UMR 7501, Université de Strasbourg et CNRS \\
7 rue René Descartes \\
67000 Strasbourg, France\\
opshtein@math.unistra.fr
\end{document}

%% file: inflation.pstex_t
\begin{picture}(0,0)%
\epsfig{file=inflation.pstex}%
\end{picture}%
\setlength{\unitlength}{1119sp}%
\begingroup\makeatletter\ifx\SetFigFont\undefined%
\gdef\SetFigFont#1#2#3#4#5{%
  \reset@font\fontsize{#1}{#2pt}%
  \fontfamily{#3}\fontseries{#4}\fontshape{#5}%
  \selectfont}%
\fi\endgroup%
\begin{picture}(8230,8226)(2096,-7954)
\put(7661,-4681){\makebox(0,0)[lb]{\smash{{\SetFigFont{5}{6.0}{\rmdefault}{\mddefault}{\updefault}{\color[rgb]{0,0,0}(II)}%
}}}}
\put(5911,-2041){\makebox(0,0)[lb]{\smash{{\SetFigFont{5}{6.0}{\rmdefault}{\mddefault}{\updefault}{\color[rgb]{0,0,0}$B(a)$}%
}}}}
\put(5376,-4071){\makebox(0,0)[lb]{\smash{{\SetFigFont{5}{6.0}{\rmdefault}{\mddefault}{\updefault}{\color[rgb]{0,0,0}$X_3(a)$}%
}}}}
\put(7141,-3181){\makebox(0,0)[lb]{\smash{{\SetFigFont{5}{6.0}{\rmdefault}{\mddefault}{\updefault}{\color[rgb]{0,0,0}$X_\beta$}%
}}}}
\put(9016,-1861){\makebox(0,0)[lb]{\smash{{\SetFigFont{5}{6.0}{\rmdefault}{\mddefault}{\updefault}{\color[rgb]{0,0,0}$\Sigma$}%
}}}}
\put(5356,-391){\makebox(0,0)[lb]{\smash{{\SetFigFont{5}{6.0}{\rmdefault}{\mddefault}{\updefault}{\color[rgb]{0,0,0}$\om_{|M\priv \Sigma}=d\beta$}%
}}}}
\put(8191,-4936){\makebox(0,0)[lb]{\smash{{\SetFigFont{5}{6.0}{\rmdefault}{\mddefault}{\updefault}{\color[rgb]{0,0,0}(III)}%
}}}}
\put(7136,-4446){\makebox(0,0)[lb]{\smash{{\SetFigFont{5}{6.0}{\rmdefault}{\mddefault}{\updefault}{\color[rgb]{0,0,0}(I)}%
}}}}
\put(5341,-61){\makebox(0,0)[lb]{\smash{{\SetFigFont{5}{6.0}{\rmdefault}{\mddefault}{\updefault}{\color[rgb]{0,0,0}$(M,\om)$}%
}}}}
\end{picture}%

%% file: fatcross2.pstex_t
\begin{picture}(0,0)%
\epsfig{file=fatcross2.pstex}%
\end{picture}%
\setlength{\unitlength}{1202sp}%
\begingroup\makeatletter\ifx\SetFigFont\undefined%
\gdef\SetFigFont#1#2#3#4#5{%
  \reset@font\fontsize{#1}{#2pt}%
  \fontfamily{#3}\fontseries{#4}\fontshape{#5}%
  \selectfont}%
\fi\endgroup%
\begin{picture}(10976,3582)(99,-5168)
\put(2129,-4761){\makebox(0,0)[lb]{\smash{{\SetFigFont{5}{6.0}{\rmdefault}{\mddefault}{\updefault}{\color[rgb]{0,0,0}$B(a_i)$}%
}}}}
\put(7316,-3290){\makebox(0,0)[lb]{\smash{{\SetFigFont{5}{6.0}{\rmdefault}{\mddefault}{\updefault}{\color[rgb]{0,0,0}$X_i$}%
}}}}
\put(4658,-2947){\makebox(0,0)[lb]{\smash{{\SetFigFont{5}{6.0}{\rmdefault}{\mddefault}{\updefault}{\color[rgb]{0,0,0}$\phi_i$}%
}}}}
\put(7426,-2221){\makebox(0,0)[lb]{\smash{{\SetFigFont{5}{6.0}{\rmdefault}{\mddefault}{\updefault}{\color[rgb]{0,0,0}$M$}%
}}}}
\put(2701,-3841){\makebox(0,0)[lb]{\smash{{\SetFigFont{5}{6.0}{\rmdefault}{\mddefault}{\updefault}{\color[rgb]{0,0,0}$\theta'_i$}%
}}}}
\put(1756,-3976){\makebox(0,0)[lb]{\smash{{\SetFigFont{5}{6.0}{\rmdefault}{\mddefault}{\updefault}{\color[rgb]{0,0,0}$B(\eps)$}%
}}}}
\put(8866,-3886){\makebox(0,0)[lb]{\smash{{\SetFigFont{5}{6.0}{\rmdefault}{\mddefault}{\updefault}{\color[rgb]{0,0,0}$B^\eps_i$}%
}}}}
\put(9811,-3796){\makebox(0,0)[lb]{\smash{{\SetFigFont{5}{6.0}{\rmdefault}{\mddefault}{\updefault}{\color[rgb]{0,0,0}$\theta_i$}%
}}}}
\put(1936,-2536){\makebox(0,0)[lb]{\smash{{\SetFigFont{5}{6.0}{\rmdefault}{\mddefault}{\updefault}{\color[rgb]{0,0,0}$\theta'_i$}%
}}}}
\put(9046,-2311){\makebox(0,0)[lb]{\smash{{\SetFigFont{5}{6.0}{\rmdefault}{\mddefault}{\updefault}{\color[rgb]{0,0,0}$\theta_i$}%
}}}}
\end{picture}%